\documentclass[10pt]{amsart}

\usepackage{amssymb}
\usepackage{amsthm}
\usepackage{amscd}
\usepackage{amsmath}
\usepackage{enumerate}
\usepackage{delarray}
\usepackage{hhline}

\title[Newton series and extended derivation relations for MLV's]{Newton series and extended derivation relations for multiple $L$-values}
\author{Gaku Kawashima and Tatsushi Tanaka}
\address{Graduate School of Mathematics, Nagoya University \endgraf
Chikusa-ku, Nagoya, 464-8602, Japan}
\email{m02009c@math.nagoya-u.ac.jp}
\address{Graduate School of Mathematics, Kyushu University \endgraf
Fukuoka, 812-8581, Japan}
\email{t.tanaka@math.kyushu-u.ac.jp}

\theoremstyle{definition}
\newtheorem{thm}{Theorem}[section]
\newtheorem{defn}[thm]{Definition}

\newtheorem{rmk}[thm]{Remark}
\newtheorem{cor}[thm]{Corollary}
\newtheorem{lem}[thm]{Lemma}
\newtheorem{prop}[thm]{Proposition}

\DeclareFontEncoding{OT2}{}{}
\DeclareFontSubstitution{OT2}{cmr}{m}{n}

\DeclareFontFamily{OT2}{cmr}{\hyphenchar\font45 }
\DeclareFontShape{OT2}{cmr}{m}{n}{%
   <5><6><7><8><9>gen*wncyr%
   <10><10.95><12><14.4><17.28><20.74><24.88>wncyr10}{}
\DeclareFontShape{OT2}{cmr}{b}{n}{%
   <5><6><7><8><9>gen*wncyb%
   <10><10.95><12><14.4><17.28><20.74><24.88>wncyb10}{}

\DeclareMathAlphabet{\mathcyr}{OT2}{cmr}{m}{n}
\DeclareMathAlphabet{\mathcyb}{OT2}{cmr}{b}{n}
\SetMathAlphabet{\mathcyr}{bold}{OT2}{cmr}{b}{n}

\begin{document}

\maketitle

\begin{abstract}

We investigate Newton series for truncated multiple $L$-values and thereby obtain a class of relations for multiple $L$-values. In addition, we give a formulation and a proof of extended derivation relations for multiple $L$-values. 
\end{abstract}

\tableofcontents


\section{Introduction}\label{sec1}

\noindent Let $r$ and $n$ be positive integers. For $r_1,\ldots,r_n\in\mathbb{Z}/r\mathbb{Z}$ and $k_1,\ldots,k_n\in\mathbb{N}$, two types of multiple $L$-values (MLV's), namely, shuffle type $(\mathcyr{sh})$ and harmonic type $(\ast)$, were defined in \cite{AK} by
$$ \displaystyle L^{\mathcyb{sh}}(\mathbf{k};\mathbf{r})=\lim_{m\to\infty}\sum_{m>m_1>\cdots >m_n>0}\frac{\zeta^{r_1(m_1-m_2)}\cdots \zeta^{r_{n-1}(m_{n-1}-m_n)}\zeta^{r_nm_n}}{m_1^{k_1}\cdots m_n^{k_n}} $$
and
$$ \displaystyle L^{\ast}(\mathbf{k};\mathbf{r})=\lim_{m\to\infty}\sum_{m>m_1>\cdots >m_n>0}\frac{\zeta^{r_1m_1}\cdots \zeta^{r_nm_n}}{m_1^{k_1}\cdots m_n^{k_n}}, $$
where $(\mathbf{k};\mathbf{r})$ denote the index set $(k_1,\ldots,k_n;r_1,\ldots,r_n)$ and $\zeta=\exp(2\pi i/r)$. We also define two types of non-strict MLV's by
$$ \displaystyle \overline{L}^{\mathcyb{sh}}(\mathbf{k};\mathbf{r})=\lim_{m\to\infty}\sum_{m\ge m_1\ge \cdots \ge m_n>0}\frac{\zeta^{r_1(m_1-m_2)}\cdots \zeta^{r_{n-1}(m_{n-1}-m_n)}\zeta^{r_nm_n}}{m_1^{k_1}\cdots m_n^{k_n}} $$
and
$$ \displaystyle \overline{L}^{\ast}(\mathbf{k};\mathbf{r})=\lim_{m\to\infty}\sum_{m\ge m_1\ge \cdots \ge m_n>0}\frac{\zeta^{r_1m_1}\cdots \zeta^{r_nm_n}}{m_1^{k_1}\cdots m_n^{k_n}}. $$
Note that both types of non-strict MLV's converge if $(k_1,r_1)\neq (1,0)$ (see \cite{AK} for example). 

The bound of the dimension of the `MLV-space', which is a $\mathbb{Q}$-vector space generated by MLV's, has been studied by P. Deligne and A. Goncharov in \cite{DG,G}. According to their result, there are several linear relations over $\mathbb{Q}$ among the MLV's. In \cite{AK}, a large class of relations for MLV's called extended double shuffle relations (EDSR), which contains the derivation relation, was introduced. In the present paper, using the theory of Newton series, we obtain a new class of relations for MLV's containing the extended derivation relation. 

Now let $s_1,\ldots,s_n\in\mathbb{C}^{\times}$. We define two types of `truncated' MLV's for an index set $(\mathbf{k};\mathbf{s})=(k_1,\ldots,k_n;s_1,\ldots,s_n)$ by
$$ \displaystyle \mathtt{S}^{\mathcyb{sh}}_{(\mathbf{k};\mathbf{s})}(m)=\sum_{m\ge m_1\ge \cdots \ge m_n\ge 0}\frac{s_1^{m_1-m_2}\cdots s_{n-1}^{m_{n-1}-m_n}s_n^{m_n+1}}{(m_1+1)^{k_1}\cdots (m_n+1)^{k_n}} $$
and
$$ \displaystyle \mathtt{S}^{\ast}_{(\mathbf{k};\mathbf{s})}(m)=\sum_{m\ge m_1\ge \cdots \ge m_n\ge 0}\frac{s_1^{m_1+1}\cdots s_n^{m_n+1}}{(m_1+1)^{k_1}\cdots (m_n+1)^{k_n}}, $$
which can be viewed as sequences in terms of $m\in\mathbb{Z}_{\ge 0}$. Clearly, non-strict MLV's are limits of these truncated MLV's as $m\to\infty$. More precisely, 
$$ \displaystyle \overline{L}^{\sharp}(\mathbf{k};\mathbf{r})=\lim_{m\to\infty}\mathtt{S}^{\sharp}_{(\mathbf{k};\mathbf{s})}(m) $$ 
with $s_i=\zeta^{r_i}$ $(i=1,2,\ldots ,n)$ and $\sharp =\mathcyr{sh}$ (shuffle) or $\ast$ (harmonic). 

We prove the inversion sequence of a truncated MLV is a linear combination of truncated MLV's in $\S \ref{sec2}$. The generalized Landen connection formula of the multiple polylogarithm function plays an important role in the proof, which can be considered as a generalization of the proof of the identity
\begin{equation}
\displaystyle \sum_{i=1}^m(-1)^i\binom{m}{i}\frac{1}{i}=-\sum_{i=1}^m\frac{1}{i} \label{eq1}
\end{equation}
studied by L. Euler in \cite{E}. 

The inversion sequences of truncated MLV's are required when Newton series of truncated MLV's are discussed. We construct Newton series of truncated MLV's, find several properties of these series, in particular, a functional equation, and obtain a class of relations for MLV's in $\S \ref{sec3}$. 

In $\S \ref{sec4}$, certain extensions of the derivation relation for MLV's are formulated. Several advantageous properties of `extended derivation operators' $\widehat{\partial}_n^{(c)}$ and $\partial_n^{(c)}$ are introduced and then a proof of the extended derivation relation for MLV's is established. 

Proofs of some lemmata are presented in $\S \ref{sec5}$. Data of computation using Risa/Asir, an open source general computer algebra system, are given in Appendix at the bottom of the paper. 


\section{MPL}\label{sec2}

In this section, we prove the generalized Landen connection formula for multiple polylogarithms (MPL's) and find the inversion sequence of a truncated MLV. To describe these or other facts precisely, it is convenient to use the algebraic setup on the non-commutative algebra $\mathcal{A}_{\mathbb{C}^{\times}}:=\mathbb{Q}\langle x,y_s|s\in\mathbb{C}^{\times}\rangle$ in infinitely many indeterminates $x,y_s(s\in\mathbb{C}^{\times})$. We introduce a number of operators on $\mathcal{A}_{\mathbb{C}^{\times}}$ to obtain the inversion sequence of a truncated MLV. 

\subsection{Differential Formula}\label{subsec2-1}

\noindent Let $(\mathbf{k};\mathbf{s})$ be the index set $(k_1,\ldots,k_n;s_1,\ldots,s_n)$ $(k_i\in\mathbb{N},s_i\in\mathbb{C})$. Here, $k_1+\cdots +k_n$ is the weight, and $n$ is the depth. 

Next, we introduce the strict and non-strict MPL's (of $\mathcyr{sh}$-type) for an index set $(\mathbf{k};\mathbf{s})$:
$$ \displaystyle \mathrm{Li}^{\mathcyb{sh}}_{(\mathbf{k};\mathbf{s})}(z)=\sum_{m_1>\cdots >m_n>0}\frac{s_1^{m_1-m_2}\cdots s_{n-1}^{m_{n-1}-m_n}s_n^{m_n}}{m_1^{k_1}\cdots m_n^{k_n}}z^{m_1} $$
and
$$ \displaystyle \overline{\mathrm{Li}}^{\mathcyb{sh}}_{(\mathbf{k};\mathbf{s})}(z)=\sum_{m_1\ge \cdots \ge m_n>0}\frac{s_1^{m_1-m_2}\cdots s_{n-1}^{m_{n-1}-m_n}s_n^{m_n}}{m_1^{k_1}\cdots m_n^{k_n}}z^{m_1} $$
where $z$ is a complex variable. These MPL's (with $s_1\neq 0$) converge absolutely if $|z|<|1/s_1|$. For the index set $(1;s)$ ($s\in\mathbb{C}$), we have
$$ \displaystyle \mathrm{Li}^{\mathcyb{sh}}_{(1;s)}(z)=-\log (1-sz) $$
and the identity
\begin{equation}
\displaystyle \mathrm{Li}^{\mathcyb{sh}}_{(1;s)}\biggl(\frac{z}{1-z}\biggl)=\mathrm{Li}^{\mathcyb{sh}}_{(1;1-s)}(z)-\mathrm{Li}^{\mathcyb{sh}}_{(1;1)}(z) \label{eq2}
\end{equation}
which are fundamental properties of the logarithm function. 

\begin{lem}[Differential formula]\label{lem2-1}
{\it For $|z|<|1/s_1|$ and an index set $(\mathbf{k};\mathbf{s})$ consisting of $k_1,\ldots,k_n\in\mathbb{N}$ and $s_1,\ldots,s_n\in\mathbb{C}^{\times}$, we have
\[
\displaystyle\frac{d}{dz}\mathrm{Li}^{\mathcyb{sh}}_{(\mathbf{k};\mathbf{s})}(z)=
\left\{
\begin{array}{ll}
\displaystyle\frac{1}{z}\mathrm{Li}^{\mathcyb{sh}}_{(k_1-1,k_2,\ldots,k_n;\mathbf{s})}(z) & k_1>1, \\
\displaystyle\frac{s_1}{1-s_1z}\mathrm{Li}^{\mathcyb{sh}}_{(k_2,\ldots,k_n;s_2,\ldots,s_n)}(z) & k_1=1,~n>1, \\
\displaystyle\frac{s_1}{1-s_1z} & k_1=1,~n=1
\end{array}
\right.
\]
and
\begin{eqnarray*}
\displaystyle &{}& \frac{d}{dz}\mathrm{Li}^{\mathcyb{sh}}_{(\mathbf{k};\mathbf{s})}\biggl(\frac{z}{1-z}\biggl) \\
&{}& =\left\{
\begin{array}{ll}
\displaystyle\biggl(\frac{1}{z}+\frac{1}{1-z}\biggl)\mathrm{Li}^{\mathcyb{sh}}_{(k_1-1,k_2,\ldots,k_n;\mathbf{s})}\biggl(\frac{z}{1-z}\biggl) & k_1>1, \\ 
\displaystyle\biggl(\frac{1-s_1}{1-(1-s_1)z}-\frac{1}{1-z}\biggl)\mathrm{Li}^{\mathcyb{sh}}_{(k_2,\ldots,k_n;s_2,\ldots,s_n)}\biggl(\frac{z}{1-z}\biggl) & k_1=1, n>1, \\
\displaystyle\frac{1-s_1}{1-(1-s_1)z}-\frac{1}{1-z} & k_1=1, n=1. 
\end{array}
\right.
\end{eqnarray*}
}
\end{lem}

\begin{proof}
The proof is simple for $k_1>1$ and for $k_1=1$, $n=1$. If $k_1=1$, $n>1$, 

$\displaystyle \frac{d}{dz}\mathrm{Li}^{\mathcyb{sh}}_{(1,k_2,\ldots ,k_n;s_1,\ldots ,s_n)}(z) $
\begin{eqnarray*}
\displaystyle &=& \sum_{m_1>\cdots >m_n>0}\frac{s_1^{m_1-m_2}\cdots s_{n-1}^{m_{n-1}-m_n}s_n^{m_n}}{m_2^{k_2}\cdots m_n^{k_n}}z^{m_1-1} \\
\displaystyle &=& \sum_{m_2>\cdots >m_n>0}\Biggl(\sum_{m_1=m_2+1}^{\infty}s_1^{m_1}z^{m_1-1}\Biggl)\frac{s_1^{-m_2}s_2^{m_2-m_3}\cdots s_{n-1}^{m_{n-1}-m_n}s_n^{m_n}}{m_2^{k_2}\cdots m_n^{k_n}} \\
\displaystyle &=& \sum_{m_2>\cdots >m_n>0}\frac{s_1^{m_2+1}z^{m_2}}{1-s_1z}\frac{s_1^{-m_2}s_2^{m_2-m_3}\cdots s_{n-1}^{m_{n-1}-m_n}s_n^{m_n}}{m_2^{k_2}\cdots m_n^{k_n}} \\
\displaystyle &=& \frac{s_1}{1-s_1z}\mathrm{Li}^{\mathcyb{sh}}_{(k_2,\ldots,k_n;s_2,\ldots,s_n)}(z).
\end{eqnarray*}
Using $\frac{d}{dz}(\frac{z}{z-1})=-\frac{1}{(z-1)^2}$, we obtain the second formula. 
\end{proof}

\subsection{Algebraic Setup}\label{subsec2-2}

\noindent In order to provide precise descriptions, we use an algebraic setup that is similar to that of Arakawa-Kaneko \cite{AK}. Let $\Lambda$ be a group. We denote by $\mathcal{A}_{\Lambda}$ the non-commutative algebra $\mathbb{Q}\langle x,y_s|s\in\Lambda\rangle$ in indeterminates $x,y_s(s\in\Lambda)$. The subalgebras $\mathcal{A}_{\Lambda}^1$ and $\mathcal{A}_{\Lambda}^0$ are given by
$$ \displaystyle \mathcal{A}_{\Lambda}\supset\mathcal{A}_{\Lambda}^1:=\mathbb{Q}+\sum_{s\in\Lambda}\mathcal{A}_{\Lambda}y_s\supset\mathcal{A}_{\Lambda}^0:=\mathbb{Q}+\sum_{s\in\Lambda}x\mathcal{A}_{\Lambda}y_s+\sum_{s,t\in\Lambda, t\neq 1}y_t\mathcal{A}_{\Lambda}y_s. $$
We often view an index set $(\mathbf{k};\mathbf{s})=(k_1,\ldots,k_n;s_1,\ldots,s_n)$ $(k_i\in\mathbb{N},s_i\in\Lambda)$ as a word $x^{k_1-1}y_{s_1}\cdots x^{k_n-1}y_{s_n}\in\mathcal{A}_{\Lambda}^1$ and vice versa. A word $w$ is called admissible if $w\in\mathcal{A}_{\Lambda}^0$, and an index set $(\mathbf{k};\mathbf{s})$ $(s\in\Lambda)$ is called admissible if the corresponding word is admissible. Here, the number of $x$ and $y_s$ ($y_s$) of a word $w$ is its weight (depth), which corresponds with the naming convention of the corresponding index set. 

The MPL-evaluation maps $\mathrm{Li}^{\mathcyb{sh}}_{\bullet}(z), \overline{\mathrm{Li}}^{\mathcyb{sh}}_{\bullet}(z):\mathcal{A}_{\Lambda}^1\to\mathbb{C}[[z]]$ are defined by $\mathbb{Q}$-linearity and
$$ \mathrm{Li}^{\mathcyb{sh}}_{w}(z)=\mathrm{Li}^{\mathcyb{sh}}_{(\mathbf{k};\mathbf{s})}(z),\mathrm{Li}^{\mathcyb{sh}}_{1}(z)=1 $$
and
$$ \overline{\mathrm{Li}}^{\mathcyb{sh}}_{w}(z)=\overline{\mathrm{Li}}^{\mathcyb{sh}}_{(\mathbf{k};\mathbf{s})}(z),\overline{\mathrm{Li}}^{\mathcyb{sh}}_{1}(z)=1, $$
respectively, where $w=x^{k_1-1}y_{s_1}\cdots x^{k_n-1}y_{s_n}$.

As complex-valued sequences from $\mathbb{Z}_{\ge 0}$ to $\mathbb{C}$, two types of truncated MLV's $\mathtt{S}^{\mathcyb{sh}}_{(\mathbf{k};\mathbf{s})}$ and $\mathtt{S}^{\ast}_{(\mathbf{k};\mathbf{s})}$ were defined in $\S \ref{sec1}$. Here, we also define alternative truncated MLV's $\mathtt{s}^{\mathcyb{sh}}_{(\mathbf{k};\mathbf{s})}$ and $\mathtt{s}^{\ast}_{(\mathbf{k};\mathbf{s})}$ for the index set $(\mathbf{k};\mathbf{s})$ with $k_i\in\mathbb{N}$ and $s_i\in\mathbb{C}^{\times}$ by
$$ \displaystyle \mathtt{s}^{\mathcyb{sh}}_{(\mathbf{k};\mathbf{s})}(m)=\sum_{m=m_1\ge \cdots \ge m_n\ge 0}\frac{s_1^{m_1-m_2}\cdots s_{n-1}^{m_{n-1}-m_n}s_n^{m_n+1}}{(m_1+1)^{k_1}\cdots (m_n+1)^{k_n}} $$
and
$$ \displaystyle \mathtt{s}^{\ast}_{(\mathbf{k};\mathbf{s})}(m)=\sum_{m=m_1\ge \cdots \ge m_n\ge 0}\frac{s_1^{m_1+1}\cdots s_n^{m_n+1}}{(m_1+1)^{k_1}\cdots (m_n+1)^{k_n}}. $$
The evaluation maps of truncated MLV's $\mathtt{S}^{\sharp}_{\bullet},\mathtt{s}^{\sharp}_{\bullet}:\mathcal{A}_{\mathbb{C}^{\times}}^1\to\mathbb{C}$ are also defined by $\mathbb{Q}$-linearity and the mappings as well as the MPL-evaluation maps, that is 
$$ \mathtt{S}_w^{\sharp}(m)=\mathtt{S}_{(\mathbf{k};\mathbf{s})}^{\sharp}(m),\mathtt{S}_1^{\sharp}(m)=1 $$
and
$$ \mathtt{s}_w^{\sharp}(m)=\mathtt{s}_{(\mathbf{k};\mathbf{s})}^{\sharp}(m),\mathtt{s}_1^{\sharp}(m)=1 $$
where $w=x^{k_1-1}y_{s_1}\cdots x^{k_n-1}y_{s_n}$ and $\sharp=\mathcyr{sh}$ or $\ast$.

In addition, we introduce three operators on the space of complex-valued sequences $\mathbb{Z}_{\ge 0}\to\mathbb{C}$: the partial-sum operator $\Sigma$, its inverse operator $\Sigma^{-1}$, and the inversion operator $\nabla$, given by 
$$ \displaystyle (\Sigma a)(m)=\sum_{i=0}^ma(i), $$
$$ \displaystyle (\Sigma^{-1}a)(m)=\left\{
\begin{array}{ll}
a(0) & m=0, \\
a(m)-a(m-1) & m>0
\end{array}
\right. $$
and
$$ \displaystyle (\nabla a)(m)=\sum_{i=0}^m(-1)^i\binom{m}{i}a(i). $$
Note that $\Sigma\Sigma^{-1}=\Sigma^{-1}\Sigma=\mathrm{id}$, $\nabla^2=\mathrm{id}$ and $(\Sigma\nabla)^2=\mathrm{id}$. 

Under these notations, we find that the identity
\begin{equation}
\mathtt{S}^{\sharp}_w=\Sigma \mathtt{s}^{\sharp}_w \label{eq3}
\end{equation}
holds. We study an explicit form of the sequence $\nabla \mathtt{s}^{\mathcyb{sh}}_w$ in the next subsection.

\subsection{Generalized Landen Connection Formula for MPL's and Inversion Sequences of Truncated MLV's}\label{subsec2-3}

\noindent Let $\varphi$ and $\iota$ be automorphisms on $\mathcal{A}_{\mathbb{C}^{\times}}$ given by
$$ \varphi(x)=x+y_1,\varphi(y_s)=\delta(s)y_{s}-y_1 $$
and
$$ \iota(x)=x,\iota(y_s)=\delta(s)y_{1-s}-(1-\delta(s))y_1, $$
where $\delta$ is defined on $\mathbb{C}$ by
$$ \delta(s)=
\left\{
\begin{array}{ll}
0 & s=0,1, \\
1 & \mathrm{otherwise}.
\end{array}
\right. $$
Note that $\varphi^2=\mathrm{id}$, $\iota^2=\mathrm{id}$ and $\varphi\iota=\iota\varphi$. 

\begin{thm}[Generalized Landen Connection Formula]\label{thm2-2}
{\it For any word $w\in\mathcal{A}_{\mathbb{C}^{\times}}^1$, there exists $\varepsilon >0$ such that the identity
$$ \displaystyle \mathrm{Li}^{\mathcyb{sh}}_w(z)=\mathrm{Li}^{\mathcyb{sh}}_{\varphi\iota(w)}\biggl(\frac{z}{z-1}\biggl) $$
holds in the open disc $|z|<\varepsilon$. }
\end{thm}

\begin{proof}
Such $\varepsilon >0$ can be taken as the minimal number of radii of convergence of appearing MPL's. The identity itself is proven by the induction on the weight of the word $w$. If $w=y_s$, it is nothing but the identity (\ref{eq2}). Since $\varphi$ and $\iota$ are automorphisms, we have

\vspace{10pt}

$\displaystyle \mathrm{Li}^{\mathcyb{sh}}_{\varphi\iota(x^{k_1-1}y_{s_1}\cdots x^{k_n-1}y_{s_n})}\Bigl(\frac{z}{z-1}\Bigr)$
\begin{eqnarray*}
\displaystyle &=& 
\left\{
\begin{array}{ll}
\displaystyle \mathrm{Li}^{\mathcyb{sh}}_{(x+y_1)\varphi\iota(x^{k_1-2}y_{s_1}x^{k_2-1}y_{s_2}\cdots x^{k_n-1}y_{s_n})}\Bigl(\frac{z}{z-1}\Bigr) & k_1>1, \\
\displaystyle \mathrm{Li}^{\mathcyb{sh}}_{(\delta(s_1)y_{1-s_1}-y_1)\varphi\iota(x^{k_2-1}y_{s_2}\cdots x^{k_n-1}y_{s_n})}\Bigl(\frac{z}{z-1}\Bigr) & k_1=1.
\end{array}
\right. \\
\end{eqnarray*}
By Lemma \ref{lem2-1} and the induction hypothesis, we have

$\displaystyle \frac{d}{dz}\mathrm{Li}^{\mathcyb{sh}}_{\varphi\iota(x^{k_1-1}y_{s_1}\cdots x^{k_n-1}y_{s_n})}\Bigl(\frac{z}{z-1}\Bigr)$
\begin{eqnarray*}
\displaystyle &=& 
\left\{
\begin{array}{ll}
\displaystyle \frac{1}{z}\mathrm{Li}^{\mathcyb{sh}}_{\varphi\iota(x^{k_1-2}y_{s_1}x^{k_2-1}y_{s_2}\cdots x^{k_n-1}y_{s_n})}\Bigl(\frac{z}{z-1}\Bigr) & k_1>1 \\
\displaystyle \frac{s_1}{1-s_1z}\mathrm{Li}^{\mathcyb{sh}}_{\varphi\iota(x^{k_2-1}y_{s_2}\cdots x^{k_n-1}y_{s_n})}\Bigl(\frac{z}{z-1}\Bigr) & k_1=1,n>1 \\
\displaystyle \frac{s_1}{1-s_1z} & k_1=1,n=1
\end{array}
\right. \\
&=& 
\left\{
\begin{array}{ll}
\displaystyle \frac{1}{z}\mathrm{Li}^{\mathcyb{sh}}_{x^{k_1-2}y_{s_1}x^{k_2-1}y_{s_2}\cdots x^{k_n-1}y_{s_n}}(z) & k_1>1 \\
\displaystyle \frac{s_1}{1-s_1z}\mathrm{Li}^{\mathcyb{sh}}_{x^{k_2-1}y_{s_2}\cdots x^{k_n-1}y_{s_n}}(z) & k_1=1,n>1 \\
\displaystyle \frac{s_1}{1-s_1z} & k_1=1,n=1
\end{array}
\right. \\
&=& \frac{d}{dz}\mathrm{Li}^{\mathcyb{sh}}_{\varphi\iota(x^{k_1-1}y_{s_1}\cdots x^{k_n-1}y_{s_n})}(z). 
\end{eqnarray*}
Integrate both sides from $0$ to $z$, and we obtain Theorem.
\end{proof}

Let $\alpha$ and $\gamma$ be automorphisms on $\mathcal{A}_{\mathbb{C}^{\times}}$ and let $R_w(w\in\mathcal{A}_{\mathbb{C}^{\times}})$ be a $\mathbb{Q}$-linear operator on $\mathcal{A}_{\mathbb{C}^{\times}}$ given by
$$ \alpha(x)=y_1,\alpha(y_s)=(1-\delta(s))x+\delta(s)y_s, $$
$$ \gamma(x)=x,\gamma(y_s)=x+y_s $$
and
$$ R_w(w^{\prime})=w^{\prime}w~(w^{\prime}\in\mathcal{A}_{\mathbb{C}^{\times}}). $$
We define $\mathbb{Q}$-linear operators $\star$ and $d_{\mathcyb{sh}}$ on $\mathcal{A}^1_{\mathbb{C}^{\times}}$, respectively, by
$$ \star(wy_s)=\alpha\iota(w)(y_1-\delta(s)y_{1-s}) $$
and
$$ d_{\mathcyb{sh}}(wy_s)=\gamma(w)y_s, $$
where $w\in\mathcal{A}_{\mathbb{C}^{\times}}$. Note that $\alpha^2=\mathrm{id}$, $\alpha\iota=\iota\alpha$ and $\star^2=\mathrm{id}$. We also find that 
\begin{equation}
\overline{\mathrm{Li}}^{\mathcyb{sh}}_w(z)=\mathrm{Li}^{\mathcyb{sh}}_{d_{\mathcyb{sh}}(w)}(z). \label{eq4}
\end{equation}
\begin{lem}\label{lem2-3}
{\it $\varphi\iota d_{\mathcyb{sh}}=-d_{\mathcyb{sh}}\star.$ }
\end{lem}

\begin{proof}
For $w\in\mathcal{A}_{\mathbb{C}^{\times}}$, we have
$$ \varphi\iota d_{\mathcyb{sh}}(wy_s)=\varphi\iota(\gamma(w)y_s)=\varphi\iota\gamma(w)(\delta(s)y_{1-s}-y_1) $$
and
$$ -d_{\mathcyb{sh}}\star(wy_s)=-d_{\mathcyb{sh}}(\alpha(w)(y_1-\delta(s)y_{1-s}))=\gamma\alpha(w)(\delta(s)y_{1-s}-y_1). $$
It is a simple task to check the identity $\varphi\iota\gamma=\gamma\alpha$. 
\end{proof}

Using this lemma, we interpret the generalized Landen connection formula for the strict MPL in Theorem \ref{thm2-2} as that for the non-strict MPL. 
\begin{cor}\label{cor2-4}
{\it For any word $w\in\mathcal{A}_{\mathbb{C}^{\times}}^1$, there exists $\varepsilon >0$ such that the identity
$$ \displaystyle \overline{\mathrm{Li}}^{\mathcyb{sh}}_w(z)=-\overline{\mathrm{Li}}^{\mathcyb{sh}}_{\star(w)}\biggl(\frac{z}{z-1}\biggl) $$
holds in the open disc $|z|<\varepsilon$. }
\end{cor}

\begin{proof}
By equation (\ref{eq4}), Theorem \ref{thm2-2}, and Lemma \ref{lem2-3}, we have
$$ \displaystyle \mathrm{LHS}=\mathrm{Li}^{\mathcyb{sh}}_{d_{\mathcyb{sh}}(w)}(z)=\mathrm{Li}^{\mathcyb{sh}}_{\varphi\iota d_{\mathcyb{sh}}(w)}\Bigl(\frac{z}{z-1}\Bigr)=-\mathrm{Li}^{\mathcyb{sh}}_{d_{\mathcyb{sh}}\star(w)}\Bigl(\frac{z}{z-1}\Bigr),$$
which is made equivalent to the RHS of the identity by using equation (\ref{eq4}) again. 
\end{proof}
\begin{lem}\label{lem2-5}
{\it There exists $\varepsilon >0$ such that the identities
$$ \displaystyle \frac{1}{1-z}\overline{\mathrm{Li}}^{\mathcyb{sh}}_w(z)=\sum_{m\ge 0}(\Sigma \mathtt{s}^{\mathcyb{sh}}_w)(m)z^{m+1} $$
and
$$ \displaystyle \frac{1}{1-z}\overline{\mathrm{Li}}^{\mathcyb{sh}}_w\biggl(\frac{z}{z-1}\biggl)=-\sum_{m\ge 0}(\nabla\Sigma^{-1} \mathtt{s}^{\mathcyb{sh}}_w)(m)z^{m+1} $$
hold in the open disc $|z|<\varepsilon$. }
\end{lem}

\begin{proof}
Such $\varepsilon >0$ can be taken as the minimal number of radii of convergence of appearing MPL's. By the definition, we have
$$ \displaystyle \overline{\mathrm{Li}}^{\mathcyb{sh}}_w(z)=\sum_{m\ge 1}\mathtt{s}^{\mathcyb{sh}}_w(m-1)z^m. $$
Then, \\
$$ \displaystyle \frac{1}{1-z}\overline{\mathrm{Li}}^{\mathcyb{sh}}_w(z)=\sum_{l=1}^{\infty}\mathtt{s}^{\mathcyb{sh}}_w(l-1)\frac{z^l}{1-z}=\sum_{l=1}^{\infty}\mathtt{s}^{\mathcyb{sh}}_w(l-1)\sum_{m=l}^{\infty}z^m $$
$$ \displaystyle =\sum_{m>0}\sum_{l=1}^m \mathtt{s}^{\mathcyb{sh}}_w(l-1)z^m=\sum_{m\ge 0}(\Sigma \mathtt{s}^{\mathcyb{sh}}_w)(m)z^{m+1} $$
and 
$$ \displaystyle \frac{1}{1-z}\overline{\mathrm{Li}}^{\mathcyb{sh}}_w\Bigl(\frac{z}{z-1}\Bigr)=\sum_{l=1}^{\infty}(-1)^l \mathtt{s}^{\mathcyb{sh}}_w(l-1)\frac{z^l}{(1-z)^{l+1}}=\sum_{l=1}^{\infty}(-1)^l \mathtt{s}^{\mathcyb{sh}}_w(l-1)\sum_{m=l}^{\infty}\dbinom{m}{l}z^m $$
$$ \displaystyle =\sum_{m>0}\sum_{l=1}^m (-1)^l\dbinom{m}{l}\mathtt{s}^{\mathcyb{sh}}_w(l-1)z^m=-\sum_{m\ge 0}(\nabla\Sigma^{-1} \mathtt{s}^{\mathcyb{sh}}_w)(m)z^{m+1}.$$
\end{proof}

Using these properties of MPL's, we obtain the inversion sequences of truncated MLV's. 
\begin{thm}\label{thm2-6}
{\it For any word $w\in\mathcal{A}^{\mathcyb{sh}}_{\mathbb{C}^{\times}}$, we have 
$$ \nabla \mathtt{s}^{\mathcyb{sh}}_w=\mathtt{s}^{\mathcyb{sh}}_{\star(w)}. $$ 
}
\end{thm}

\begin{proof}
By Corollary \ref{cor2-4} and Lemma \ref{lem2-5}, we have $\Sigma \mathtt{s}^{\mathcyb{sh}}_w=\nabla\Sigma^{-1}\mathtt{s}^{\mathcyb{sh}}_{\star(w)}$. Since both $\nabla$ and $\Sigma\nabla$ are involutions, we have $\nabla \mathtt{s}^{\mathcyb{sh}}_w=\Sigma\nabla\Sigma \mathtt{s}^{\mathcyb{sh}}_w=\mathtt{s}^{\mathcyb{sh}}_{\star(w)}$. 
\end{proof}

\begin{rmk}\label{rmk2-7}
In the case of $w=y_1$, we have
\begin{equation}
\displaystyle \sum_{i=0}^m(-1)^i\binom{m}{i}\frac{1}{i+1}=\frac{1}{m+1} \label{eq5}
\end{equation}
for $m\ge 0$ because $\star(w)=y_1$ and $\mathtt{s}^{\mathcyb{sh}}_{y_1}(m)=\frac{1}{m+1}$. We can confirm that equation (\ref{eq5}) is equivalent to Euler's equation (\ref{eq1}). 
\end{rmk}


\section{Newton Series}\label{sec3}

\noindent The Newton series for a sequence $a:\mathbb{Z}_{\ge 0}\to\mathbb{C}$ is a one-variable complex function that interpolates the sequence $a$. In this section, we prove several analytic properties of the Newton series for truncated MLV's and obtain a class of relations for MLV's. 

\subsection{Order of the $l$-th Difference of Truncated MLV's}\label{subsec3-1}

\noindent For $s_1,\ldots,s_p\in\mathbb{C}$ and $m\in\mathbb{Z}_{\ge 0}$, we set\begin{equation*}
\displaystyle \mathtt{c}_{s_1,\ldots,s_p}(m)=\sum_{m=m_1\ge\cdots\ge m_p\ge 0}\frac{s_1^{m_1-m_2}\cdots s_{p-1}^{m_{p-1}-m_p}}{(m_1+1)\cdots (m_{p-1}+1)}s_p^{m_p}.
\end{equation*}
In addition, for $s_1,\ldots,s_p,t_1,\ldots,t_p\in\mathbb{C}$ and $m,l\in\mathbb{Z}_{\ge 0}$, we set

$ \displaystyle \mathtt{c}_{s_1,\ldots,s_p;t_1,\ldots,t_p}(m,l) $
$$ \displaystyle =\binom{m+l}{m}^{-1}\sum_{\begin{subarray}{c}m=m_1\ge\cdots\ge m_p\ge 0, \\ l=l_1\ge\cdots\ge l_p\ge 0 \end{subarray}}\binom{m_1-m_2+l_1-l_2}{m_1-m_2}\cdots\binom{m_{p-1}-m_p+l_{p-1}-l_p}{m_{p-1}-m_p} $$
$$ \displaystyle\qquad\qquad \times\binom{m_p+l_p}{m_p}\frac{s_1^{m_1-m_2}t_1^{l_1-l_2}\cdots s_{p-1}^{m_{p-1}-m_p}t_{p-1}^{l_{p-1}-l_p}}{(m_1+l_1+1)\cdots (m_{p-1}+l_{p-1}+1)}s_p^{m_p}t_p^{l_p}. $$
If $s_1,\ldots,s_p,t_1,\ldots,t_p\ge 0$, we see that $\mathtt{c}_{s_1,\ldots,s_p;t_1,\ldots,t_p}(m,l)\ge 0$ for any $m,l\in\mathbb{Z}_{\ge 0}$.

\begin{lem}\label{lem1}
{\it Let $p\ge 2$ and $s_1,\ldots,s_p,t_1,\ldots,t_p\in\mathbb{C}$. For any $m,l\in\mathbb{Z}_{\ge 0}$, we have}
$$ \displaystyle (m+l+1)\mathtt{c}_{s_1,\ldots,s_p;t_1,\ldots,t_p}(m,l)-s_1m\mathtt{c}_{s_1,\ldots,s_p;t_1,\ldots,t_p}(m-1,l) $$
$$ \displaystyle -t_1l\mathtt{c}_{s_1,\ldots,s_p;t_1,\ldots,t_p}(m,l-1)=\mathtt{c}_{s_2,\ldots,s_p;t_2,\ldots,t_p}(m,l). $$
\end{lem}
\begin{proof}
By definition, we have

$ \displaystyle (m+l+1)\mathtt{c}_{s_1,\ldots,s_p;t_1,\ldots,t_p}(m,l) $
$$ \displaystyle =\binom{m+l}{m}^{-1}\Biggl\{\sum_{\begin{subarray}{c}m=m_2\ge\cdots\ge m_p\ge 0, \\ l=l_2\ge\cdots\ge l_p\ge 0 \end{subarray}}\binom{m_2-m_3+l_2-l_3}{m_2-m_3}\cdots\binom{m_{p-1}-m_p+l_{p-1}-l_p}{m_{p-1}-m_p} $$
$$ \displaystyle \times\binom{m_p+l_p}{m_p}\frac{s_1^{m_1-m_2}t_1^{l_1-l_2}\cdots s_{p-1}^{m_{p-1}-m_p}t_{p-1}^{l_{p-1}-l_p}}{(m_2+l_2+1)\cdots (m_{p-1}+l_{p-1}+1)}s_p^{m_p}t_p^{l_p} $$
$$ \displaystyle +\sum_{\begin{subarray}{c}m\ge m_2\ge\cdots\ge m_p\ge 0, \\ l\ge l_2\ge\cdots\ge l_p\ge 0, \\ m>m_2~\mathrm{or}~l>l_2 \end{subarray}}\binom{m-m_2+l-l_2}{m-m_2}\binom{m_2-m_3+l_2-l_3}{m_2-m_3}\cdots\binom{m_{p-1}-m_p+l_{p-1}-l_p}{m_{p-1}-m_p} $$
$$ \displaystyle \times\binom{m_p+l_p}{m_p}\frac{s_1^{m_1-m_2}t_1^{l_1-l_2}\cdots s_{p-1}^{m_{p-1}-m_p}t_{p-1}^{l_{p-1}-l_p}}{(m_2+l_2+1)\cdots (m_{p-1}+l_{p-1}+1)}s_p^{m_p}t_p^{l_p}\Biggr\}. $$
Using the identity 
$$ \displaystyle \binom{m-m_2+l-l_2}{m-m_2}=\binom{m-m_2+l-l_2-1}{m-m_2}+\binom{m-m_2+l-l_2-1}{m-m_2-1}, $$
we have
$$ \displaystyle \mathrm{RHS}=\mathtt{c}_{s_2,\ldots,s_p,t_2,\ldots,t_p}(m,l)+\binom{m+l}{m}^{-1}\Biggl\{\sum_{\begin{subarray}{c}m\ge m_2\ge\cdots\ge m_p\ge 0, \\ l>l_2\ge\cdots\ge l_p\ge 0 \end{subarray}}\binom{m-m_2+l-l_2-1}{m-m_2}Q $$
$$ \displaystyle +\sum_{\begin{subarray}{c}m>m_2\ge\cdots\ge m_p\ge 0, \\ l\ge l_2\ge\cdots\ge l_p\ge 0 \end{subarray}}\binom{m-m_2+l-l_2-1}{m-m_2-1}Q\Biggr\}  $$
where
$$ \displaystyle Q=\binom{m_2-m_3+l_2-l_3}{m_2-m_3}\cdots\binom{m_{p-1}-m_p+l_{p-1}-l_p}{m_{p-1}-m_p}\binom{m_p+l_p}{m_p} $$
$$ \displaystyle \times\frac{t_1^{l-l_2}s_2^{m_2-m_3}t_2^{l_2-l_3}\cdots s_{p-1}^{m_{p-1}-m_p}t_{p-1}^{l_{p-1}-l_p}}{(m_2+l_2+1)\cdots (m_{p-1}+l_{p-1}+1)}s_p^{m_p}t_p^{l_p}. $$
This implies the identity of Lemma. 
\end{proof}

We introduce an operator on the space of complex-valued sequences called the difference operator $\Delta$ given by
$$ (\Delta a)(m)=a(m)-a(m+1) $$
for any sequence $a:\mathbb{Z}_{\ge 0}\to\mathbb{C}$. 
\begin{prop}\label{prop2}
{\it Let $s_1,\ldots,s_p\in\mathbb{C}$ and $m,l\in\mathbb{Z}_{\ge 0}$. We have}
$$ (\Delta^l\mathtt{c}_{s_1,\ldots,s_p})(m)=\mathtt{c}_{s_1,\ldots,s_p;1-s_1,\ldots,1-s_p}(m,l). $$
\end{prop}
\begin{proof}
Set the generating functions
$$ \displaystyle f_{s_1,\ldots,s_p}(X,Y)=\sum_{m,l=0}^{\infty}(\Delta^l\mathtt{c}_{s_1,\ldots,s_p})(m)\frac{X^mY^l}{m!l!}, $$
$$ \displaystyle g_{s_1,\ldots,s_p}(X,Y)=\sum_{m,l=0}^{\infty}\mathtt{c}_{s_1,\ldots,s_p;1-s_1,\ldots,1-s_p}(m,l)\frac{X^mY^l}{m!l!}. $$
Then, we need only show that 
$$ f_{s_1,\ldots,s_p}(X,Y)-g_{s_1,\ldots,s_p}(X,Y)=0, $$
which is equivalent to showing that 
\[
\left\{
\begin{array}{ll}
\displaystyle \mathrm{i)}~f_{s_1,\ldots,s_p}(X,0)-g_{s_1,\ldots,s_p}(X,0)=0, \\
\displaystyle \mathrm{ii)}~(\partial_X+\partial_Y-1)(f_{s_1,\ldots,s_p}(X,Y)-g_{s_1,\ldots,s_p}(X,Y))=0. 
\end{array}
\right.
\]
Identity i) is trivial. Since $(\partial_X+\partial_Y-1)f_{s_1,\ldots,s_p}(X,Y)=0$, we need only show
\begin{equation}
(\partial_X+\partial_Y-1)g_{s_1,\ldots,s_p}(X,Y)=0. \label{eq1}
\end{equation}
The proof is by induction on $p$. Since $g_{s_1}(X,Y)=e^{s_1X+(1-s_1)Y}$, the identity (\ref{eq1}) holds. Suppose that $p\ge 2$. By Lemma \ref{lem1}, we have
\begin{equation}
(X\partial_X+Y\partial_Y+1-s_1X-(1-s_1)Y)g_{s_1\ldots s_p}(X,Y)=g_{s_2\ldots s_p}(X,Y). \label{eq2}
\end{equation}
Based on the identity (\ref{eq2}), the identity
$$ [\partial_X+\partial_Y-1,X\partial_X+Y\partial_Y+1-s_1X-(1-s_1)Y]=\partial_X+\partial_Y-1, $$
and the induction hypothesis, we have

$ (X\partial_X+Y\partial_Y+2-s_1X-(1-s_1)Y)(\partial_X+\partial_Y-1)g_{s_1\ldots s_p}(X,Y) $
$$ =(\partial_X+\partial_Y-1)(X\partial_X+Y\partial_Y+1-s_1X-(1-s_1)Y)g_{s_1\ldots s_p}(X,Y)=0. $$
Since the map $X\partial_X+Y\partial_Y+2-s_1X-(1-s_1)Y:\mathbb{C}[[X,Y]]\to\mathbb{C}[[X,Y]]$ is injective, we obtain $(\partial_X+\partial_Y-1)g_{s_1,\ldots,s_p}(X,Y)=0$. 
\end{proof}

This proposition implies that, for $0\le s_i\le 1$ ($0\le i\le 1$), we have
\begin{equation}
(\Delta^l\mathtt{c}_{s_1,\ldots,s_p})(m)\ge 0 \label{EQ}
\end{equation}
for any $m,l\ge 0$. 
\begin{prop}\label{prop3}
{\it Let $l\in\mathbb{Z}_{\ge 0}$. Suppose that $|s_i|<1$ or $s_i=1$ for any $0\le i\le 1$. For any $\varepsilon >0$, we have
$$ \displaystyle (\Delta^l\mathtt{c}_{s_1,\ldots,s_p,0})(m)=\mathit{O}\biggl(\frac{1}{m^{l+1-\varepsilon}}\biggr) $$
as $m\to\infty$.}
\end{prop}
\begin{proof}
Let $t_i=1-s_i$ ($1\le i\le p$). By Proposition \ref{prop2}, 

$ \displaystyle |(\Delta^l \mathtt{c}_{s_1,\ldots,s_p,0})(m)|=|\mathtt{c}_{s_1,\ldots,s_p,0;t_1,\ldots,t_p,1}(m,l)| $
$$ \displaystyle =\Biggl|\sum_{\begin{subarray}{c}m=m_1\ge\cdots m_p\ge 0, \\ l=l_1\ge\cdots l_{p+1}\ge 0\end{subarray}}\binom{m_1-m_2+l_1-l_2}{m_1-m_2}\cdots\binom{m_{p-1}-m_p+l_{p-1}-l_p}{m_{p-1}-m_p}\binom{m_p+l_p-l_{p+1}}{m_p} $$
$$ \displaystyle \times\frac{s_1^{m_1-m_2}t_1^{l_1-l_2}\cdots s_{p-1}^{m_{p-1}-m_p}t_{p-1}^{l_{p-1}-l_p}s_p^{m_p}t_p^{l_p-l_{p+1}}}{(M_1+\cdots +M_p+L_1+\cdots +L_{p+1}+1)\cdots (M_p+L_p+L_{p+1}+1)}\Biggr| $$
$$ \displaystyle \le\sum_{\begin{subarray}{c}m=M_1+\cdots +M_p, \\ l=L_1+\cdots +L_{p+1}, \\ M_i\ge 0(1\le i\le p), \\ L_j\ge 0(1\le j\le p+1) \end{subarray}}\binom{M_1+L_1}{M_1}\cdots\binom{M_p+L_p}{M_p}\prod_{\nu=1}^p\frac{|s_{\nu}|^{M_{\nu}}|t_{\nu}|^{L_{\nu}}}{M_{\nu}+\cdots +M_p+L_{\nu}+\cdots +L_{p+1}+1}. $$
Set $I=\{1\le i\le p|s_i=1\}$ and $J=\{1\le i\le p|s_i\neq 1\}$. Then, 
$$ \displaystyle \mathrm{RHS}=\sum_{\begin{subarray}{c}m=M_1+\cdots +M_p, \\ l=L_1+\cdots +L_{p+1}, \\ M_i\ge 0(1\le i\le p), \\ L_j\ge 0(1\le j\le p+1), \\ i\in I\Rightarrow L_i=0 \end{subarray}}\Biggl\{\prod_{i\in I}\frac{1}{M_i+\cdots +M_p+L_i+\cdots +L_{p+1}+1}\Biggr\} $$
$$ \displaystyle \times\Biggl\{\prod_{i\in J}\binom{M_i+L_i}{M_i}\frac{|s_i|^{M_i}|t_i|^{L_i}}{M_i+\cdots +M_p+L_i+\cdots +L_{p+1}+1}\Biggr\}. $$
Since $|t_i|\le 1+|s_i|\le 2$ ($1\le i\le p$), 
$$ \displaystyle \mathrm{RHS}\le\sum_{\begin{subarray}{c}m=M_1+\cdots +M_p, \\ l=L_1+\cdots +L_{p+1}, \\ M_i\ge 0(1\le i\le p), \\ L_j\ge 0(1\le j\le p+1), \\ i\in I\Rightarrow L_i=0 \end{subarray}}\Biggl\{\prod_{i\in I}\frac{1}{M_i+\cdots +M_p+1}\Biggr\}\Biggl\{\prod_{i\in J}\binom{M_i+l}{M_i}\frac{|s_i|^{M_i}2^{L_i}}{M_i+\cdots +M_p+1}\Biggr\} $$
$$ \displaystyle \le 2^l\sum_{\begin{subarray}{c}m=M_1+\cdots +M_p, \\ l=L_1+\cdots +L_{p+1}, \\ M_i\ge 0(1\le i\le p), \\ L_j\ge 0(1\le j\le p+1), \\ i\in I\Rightarrow L_i=0 \end{subarray}}\Biggl\{\prod_{i\in I}\frac{1}{M_i+\cdots +M_p+1}\Biggr\}\Biggl\{\prod_{i\in J}\binom{M_i+l}{M_i}\frac{|s_i|^{M_i}}{M_i+\cdots +M_p+1}\Biggr\} $$
\begin{equation}
\displaystyle \le 2^l(l+1)^{\# J+1}\sum_{\begin{subarray}{c}m=M_1+\cdots +M_p, \\ M_i\ge 0(1\le i\le p) \end{subarray}}\Biggl\{\prod_{i\in I}\frac{1}{M_i+\cdots +M_p+1}\Biggr\}\Biggl\{\prod_{i\in J}\binom{M_i+l}{M_i}\frac{|s_i|^{M_i}}{M_i+\cdots +M_p+1}\Biggr\}. \label{id1}
\end{equation}
If $1\in I$, then
$$ \displaystyle (\ref{id1})=\frac{2^l(l+1)^{\# J+1}}{m+1}\sum_{\begin{subarray}{c}m=M_1+\cdots +M_p, \\ M_i\ge 0(1\le i\le p) \end{subarray}}\Biggl\{\prod_{i\in I\backslash\{1\}}\frac{1}{M_i+\cdots +M_p+1}\Biggr\}\Biggl\{\prod_{i\in J}\binom{M_i+l}{M_i}\frac{|s_i|^{M_i}}{M_i+\cdots +M_p+1}\Biggr\} $$
$$ \displaystyle \le \frac{2^l(l+1)^{\# J+1}}{m+1}\sum_{\begin{subarray}{c}m=M_1+\cdots +M_p, \\ M_i\ge 0(1\le i\le p) \end{subarray}}\Biggl\{\prod_{i\in I\backslash\{1\}}\frac{1}{M_i+1}\Biggr\}\Biggl\{\prod_{i\in J}\binom{M_i+l}{M_i}|s_i|^{M_i}\Biggr\} $$
$$ \displaystyle \le \frac{2^l(l+1)^{\# J+1}}{m+1}\sum_{M_2,\ldots,M_p=0}^m\Biggl\{\prod_{i\in I\backslash\{1\}}\frac{1}{M_i+1}\Biggr\}\Biggl\{\prod_{i\in J}\binom{M_i+l}{M_i}|s_i|^{M_i}\Biggr\} $$
\begin{equation}
\displaystyle \le \frac{2^l(l+1)^{\# J+1}}{m+1}\Biggl\{\prod_{i\in I\backslash\{1\}}\Biggl(\sum_{M_i=0}^m\frac{1}{M_i+1}\Biggr)\Biggr\}\Biggl\{\prod_{i\in J}\Biggl(\sum_{M_i=0}^m\binom{M_i+l}{M_i}|s_i|^{M_i}\Biggr)\Biggr\}. \label{id2}
\end{equation}
For $0\le x<1$, we have
$$ \displaystyle \sum_{M=0}^{\infty}\binom{M+l}{M}x^M=\frac{1}{l!}\sum_{M=0}^{\infty}(M+1)\cdots (M+l)x^M<\infty, $$
and hence
$$ \displaystyle (\ref{id2})=O\Biggl(\frac{(\log m)^{\# I-1}}{m}\Biggr) $$
as $m\to\infty$. If $1\in J$, then
$$ \displaystyle (\ref{id1})\le 2^l(l+1)^{\# J+1}\frac{1}{m+1}\sum_{\begin{subarray}{c}m=M_1+\cdots +M_p, \\ M_i\ge 0(1\le i\le p) \end{subarray}}\Biggl\{\prod_{i\in I}\frac{1}{M_i+1}\Biggr\}\Biggl\{\prod_{i\in J}\binom{M_i+l}{M_i}|s_i|^{M_i}\Biggr\} $$
$$ \displaystyle \le 2^l(l+1)^{\# J+1}\frac{1}{m+1}\sum_{M_1,\ldots,M_p=0}^m\Biggl\{\prod_{i\in I}\frac{1}{M_i+1}\Biggr\}\Biggl\{\prod_{i\in J}\binom{M_i+l}{M_i}|s_i|^{M_i}\Biggr\} $$
$$ \displaystyle \le 2^l(l+1)^{\# J+1}\frac{1}{m+1}\Biggl\{\prod_{i\in I}\Biggl(\sum_{M_i=0}^m\frac{1}{M_i+1}\Biggr)\Biggr\}\Biggl\{\prod_{i\in J}\Biggl(\sum_{M_i=0}^m\binom{M_i+l}{M_i}|s_i|^{M_i}\Biggr)\Biggr\} $$
$$ \displaystyle =O\Biggl(\frac{(\log m)^{\# I}}{m}\Biggr) $$
as $m\to\infty$. Thus, we conclude the Proposition. 
\end{proof}
\begin{cor}\label{cor4}
{\it Let $l\in\mathbb{Z}_{\ge 0}$. Suppose that $|s_i|<1$ or $s_i=1$ for any $0\le i\le 1$. Then, for any $\varepsilon >0$, we have
$$ \displaystyle (\Delta^l\mathtt{s}^{\mathcyb{sh}}_{(1,\ldots,1;s_1,\ldots s_p)})(m)=\mathit{O}\biggl(\frac{1}{m^{l+1-\varepsilon}}\biggr) $$
as $m\to\infty$.}
\end{cor}
\begin{proof}
The proof follows from Proposition \ref{prop3} and $\mathtt{s}^{\mathcyb{sh}}_{(1,\ldots,1;s_1,\ldots s_p)}=s_p\mathtt{c}_{s_1,\ldots,s_p,0}$. 
\end{proof}

\subsection{Basic Properties}\label{subsec3-2}

\noindent Let $a:\mathbb{Z}_{\ge 0}\to\mathbb{C}$ be a complex-valued sequence. The Newton series for the sequence $a$ is defined by
$$ \displaystyle f_a(z):=\sum_{n=0}^{\infty}(-1)^n(\nabla a)(n)\binom{z}{n}, $$
where $\binom{z}{n}=\frac{z(z-1)\cdots (z-n+1)}{n!}$ and $z$ is a complex variable. We find that $f_a(m)=a(m)$ holds for any $m\in\mathbb{Z}_{\ge 0}$. In this sense, we may denote $f_a(z)$ by $a(z)$. 

The following properties are basic in the theory of the Newton series (see \cite{Ge,I} for details). 

\begin{prop}\label{prop3-1}
{\it Let $a:\mathbb{Z}_{\ge 0}\to\mathbb{C}$ be a sequence and $z\in\mathbb{C}\backslash\mathbb{Z}_{\ge 0}$. Then, the series 
\begin{equation}
\displaystyle \sum_{n=0}^{\infty}(-1)^na(n)\binom{z}{n} \label{eq6}
\end{equation}
and the Dirichlet series
$$ \displaystyle \sum_{n=1}^{\infty}\frac{a(n)}{n^{z+1}} $$
possess one and the same abscissa of convergence and absolute convergence. }
\end{prop}
\begin{cor}\label{cor0-3-3}
{\it If $a(n)=O(1/n^{\varepsilon})$ as $n\to\infty$, the Newton series
$$ \displaystyle \sum_{n=0}^{\infty}(-1)^na(n)\binom{z}{n} $$
converges absolutely for $\mathrm{Re}(z)>-\varepsilon$.}
\end{cor}
\begin{prop}\label{prop3-2}
{\it Let $a:\mathbb{Z}_{\ge 0}\to\mathbb{C}$ be a sequence. If there exists $N\in\mathbb{Z}_{\ge 0}$ such that $a(n)=0$ for $\mathbb{Z}_{\ge 0}\ni n\ge N$, then we have $f_a(z)=0$ in the entire right half plane of convergence. }
\end{prop}
\begin{lem}\label{lem0-3-5}
{\it Denote the abscissa of convergence of Newton series
$$ \displaystyle f(z)=\sum_{n=0}^{\infty}(-1)^na(n)\binom{z}{n} $$
by $\rho$. Let $l\in\mathbb{Z}_{\ge 0}$. Then, we have
$$ \displaystyle (-1)^l\binom{z}{l}\sum_{n=0}^{\infty}(-1)^na(n)\binom{z}{n}=\sum_{n=l}^{\infty}(-1)^n\binom{n}{l}(\Delta^l a)(n-l)\binom{z}{n} $$
for $\mathrm{Re}(z)>\rho +l$.}
\end{lem}
\begin{prop}\label{prop0-3-6}
{\it Set the Newton series
$$ \displaystyle f(z)=\sum_{n=0}^{\infty}(-1)^na(n)\binom{z}{n}, g(z)=\sum_{n=0}^{\infty}(-1)^nb(n)\binom{z}{n}. $$
Denote by $\rho$ the abscissa of convergence of $f(z)$. Let $\varepsilon >0$, and suppose that the sequences $a,b:\mathbb{Z}_{\ge 0}\to\mathbb{C}$ satisfy the following conditions:
\begin{itemize}
\item[{\rm i)}] $a(m),(\Delta^l b)(m)\ge 0$ for any $m,l\in\mathbb{Z}_{\ge 0}$,
\item[{\rm ii)}] $\rho <0$,
\item[{\rm iii)}] For any $l\in\mathbb{Z}_{\ge 0}$, $(\Delta^l b)(m)=O(1/m^{l+\varepsilon})$ as $m\to\infty$.
\end{itemize}
Then, the product $f(z)g(z)$ is expressed as a Newton series that converges for $\mathrm{Re}(z)>\max\{\rho,-\varepsilon\}$.}
\end{prop}

\begin{proof}
Let $\rho^{\prime}$ denote the abscissa of convergence of $g(z)$. By Lemma \ref{lem0-3-5}, 
\begin{equation}
\displaystyle (-1)^l\binom{z}{l}g(z)=\sum_{n=l}^{\infty}(-1)^n\binom{n}{l}(\Delta^l b)(n-l)\binom{z}{n} \label{id3}
\end{equation}
for $\mathrm{Re}(z)>\rho^{\prime}+l$. According to iii), 
$$ \displaystyle \binom{n}{l}(\Delta^l b)(n-l)=O\biggl(\frac{1}{n^{\varepsilon}}\biggr)~~(n\to\infty) $$
for any $l\in\mathbb{Z}_{\ge 0}$, and hence, by Corollary \ref{cor0-3-3}, the right-hand side of (\ref{id3}) converges for $\mathrm{Re}(z)>-\varepsilon$. Since the left-hand side of (\ref{id3}) also converges for $\mathrm{Re}(z)>-\varepsilon$, the identity (\ref{id3}) holds for $\mathrm{Re}(z)>-\varepsilon$. Then, we have 
\begin{equation}
\displaystyle f(z)g(z)=\sum_{l=0}^{\infty}\sum_{n=l}^{\infty}a(l)(-1)^n\binom{n}{l}(\Delta^l b)(n-l)\binom{z}{n} \label{id4}
\end{equation}
for $\mathrm{Re}(z)>\max\{\rho,-\varepsilon \}$. Suppose that $z\in\mathbb{R}$ satisfies $\max\{\rho,-\varepsilon \}<z<0$. Then, each term of the right-hand side of (\ref{id4}) is non-negative, and hence we have
$$ \displaystyle f(z)g(z)=\sum_{n=0}^{\infty}(-1)^n\biggl\{\sum_{l=0}^n a(l)\binom{n}{l}(\Delta^l b)(n-l)\biggr\}\binom{z}{n}. $$
This shows the Proposition. 
\end{proof}

\subsection{Algebraic Preliminary}\label{subsec3-3}

\noindent Let $\Lambda$ be a group, and let $z_{k,s}$ denote $x^{k-1}y_s\in\mathcal{A}_{\Lambda}$. Note that every word $w$ with $\deg(w)>0$ in $\mathcal{A}_{\Lambda}$ can be expressed as $z_{k_1,s_1}\cdots z_{k_n,s_n}x^l$ for $k_i\in\mathbb{N},s_i\in\Lambda(1\le i\le n)$ and some $l\ge 0$. Here, we  introduce five operators: $\mathcal{I},\mathcal{I}^{-1},M_s(s\in\Lambda),L_w(w\in\mathcal{A}_{\Lambda})$, and $d_{\ast}$. The operators $\mathcal{I},\mathcal{I}^{-1},M_s(s\in\Lambda)$, and $L_w(w\in\mathcal{A}_{\Lambda})$ are $\mathbb{Q}$-linear maps on $\mathcal{A}_{\Lambda}$ defined by
$$ \mathcal{I}(z_{k_1,s_1}\cdots z_{k_n,s_n}x^l)=z_{k_1,s_1}z_{k_2,s_1s_2}\cdots z_{k_n,s_1\cdots s_n}x^l, $$
$$ \mathcal{I}^{-1}(z_{k_1,s_1}\cdots z_{k_n,s_n}x^l)=z_{k_1,s_1}z_{k_2,\frac{s_2}{s_1}}\cdots z_{k_n,\frac{s_n}{s_{n-1}}}x^l, $$
$$ M_s(z_{k_1,s_1}\cdots z_{k_n,s_n}x^l)=z_{k_1,ss_1}z_{k_2,s_2}\cdots z_{k_n,s_n}x^l $$
and
$$ L_w(w^{\prime})=ww^{\prime}~(w^{\prime}\in\mathcal{A}_{\Lambda}). $$
The operator $d_{\ast}$ is a $\mathbb{Q}$-linear map on $\mathcal{A}_{\Lambda}^1$ defined by
$$ d_{\ast}=\mathcal{I}^{-1}d_{\mathcyb{sh}}\mathcal{I} $$
where $d_{\mathcyb{sh}}$ has been given in $\S \ref{subsec2-3}$. 
\begin{lem}\label{lem3-5}
(i) $\mathcal{I}\mathcal{I}^{-1}=\mathcal{I}^{-1}\mathcal{I}=\mathrm{id}$. \\
(ii) $M_sM_t=M_{st}$ $(s,t\in\Lambda)$. \\
(iii) $\mathcal{I}L_{z_{k,s}}=L_{z_{k,s}}\mathcal{I}M_s$. \\
(iv) {\it Each of $\mathcal{I},\mathcal{I}^{-1}$. and $M_s$ commutes with $L_x$.} \\
(v) {\it $d_{\ast}$ commutes with $M_s$. }
\end{lem}

\begin{proof}
The proof of statements (i) through (iv) is simple. Let $N_s=\mathcal{I}M_s\mathcal{I}^{-1}$. We see that $N_s$ is a $\mathbb{Q}$-linear map on $\mathcal{A}_{\Lambda}$ given by
$$ N_s(z_{k_1,s_1}\cdots z_{k_n,s_n}x^l)=z_{k_1,ss_1}\cdots z_{k_n,ss_n}x^l. $$
The operators $N_s$ and $d_{\mathcyb{sh}}$ commute because
\begin{eqnarray*}
d_{\mathcyb{sh}}N_s(z_{k_1,s_1}\cdots z_{k_n,s_n}) &=& d_{\mathcyb{sh}}(z_{k_1,ss_1}\cdots z_{k_n,ss_n}) \\
&=& x^{k_1-1}(x+y_{ss_1})\cdots x^{k_{n-1}-1}(x+y_{ss_{n-1}})x^{k_n-1}y_{ss_n}
\end{eqnarray*}
and
\begin{eqnarray*}
N_sd_{\mathcyb{sh}}(z_{k_1,s_1}\cdots z_{k_n,s_n}) &=& N_f(x^{k_1-1}(x+y_{s_1})\cdots x^{k_{n-1}-1}(x+y_{s_{n-1}})x^{k_n-1}y_{s_n}) \\
&=& x^{k_1-1}(x+y_{ss_1})\cdots x^{k_{n-1}-1}(x+y_{ss_{n-1}})x^{k_n-1}y_{ss_n}.
\end{eqnarray*}
Then,
$$ d_{\ast}M_s=\mathcal{I}^{-1}d_{\mathcyb{sh}}\mathcal{I}M_s=\mathcal{I}^{-1}d_{\mathcyb{sh}}N_s\mathcal{I}=\mathcal{I}^{-1}N_sd_{\mathcyb{sh}}\mathcal{I}=M_s\mathcal{I}^{-1}d_{\mathcyb{sh}}\mathcal{I}=M_sd_{\ast} $$
and we obtain (v). 
\end{proof}

Note that the operator $\mathcal{I}$ plays a role in translation between $\mathcyr{sh}$-type and $\ast$-type. For instance, identities
\begin{equation}
\mathtt{S}^{\ast}_w(z)=\mathtt{S}^{\mathcyb{sh}}_{\mathcal{I}(w)}(z),~\mathtt{s}^{\ast}_w(z)=\mathtt{s}^{\mathcyb{sh}}_{\mathcal{I}(w)}(z) \label{eq7}
\end{equation}
hold. 

Next, we define harmonic products. 
\begin{defn}\label{defn3-6}
{\it Let $\Lambda$ be a group. We define the harmonic product $\ast:\mathcal{A}^1_{\Lambda}\times\mathcal{A}^1_{\Lambda}\to\mathcal{A}^1_{\Lambda}$ by $\mathbb{Q}$-bilinearity and
\begin{itemize}
\item[(i)] $1\ast w=w\ast 1=w$, for any $w\in\mathcal{A}^1_{\Lambda}$, 
\item[(ii)] $z_{k,s}w\ast z_{l,t}w^{\prime}=z_{k,s}(w\ast z_{l,t}w^{\prime})+z_{l,t}(z_{k,s}w\ast w^{\prime})+z_{k+l,st}(w\ast w^{\prime})$ \\
for any $k,l\ge 1$ and any words $w,w^{\prime}\in\mathcal{A}^1_{\Lambda}$. 
\end{itemize}
We define another harmonic product $\overline{\ast}:\mathcal{A}^1_{\Lambda}\times\mathcal{A}^1_{\Lambda}\to\mathcal{A}^1_{\Lambda}$ by $\mathbb{Q}$-bilinearity and
\begin{itemize}
\item[(i)] $1~\overline{\ast}~w=w~\overline{\ast}~1=w$, for any $w\in\mathcal{A}^1_{\Lambda}$, 
\item[(ii)] $z_{k,s}w~\overline{\ast}~z_{l,t}w^{\prime}=z_{k,s}(w~\overline{\ast}~z_{l,t}w^{\prime})+z_{l,t}(z_{k,s}w~\overline{\ast}~w^{\prime})-z_{k+l,st}(w~\overline{\ast}~w^{\prime})$ \\
for any $k,l\ge 1$ and any words $w,w^{\prime}\in\mathcal{A}^1_{\Lambda}$. 
\end{itemize}
In addition, we define two harmonic products $\dot{\ast},\dot{\overline{\ast}}:\mathcal{A}^1_{\Lambda}\times\mathcal{A}^1_{\Lambda}\to\mathcal{A}^1_{\Lambda}$  by $\mathbb{Q}$-bilinearity and
\begin{itemize}
\item[(i)] $1~\dot{\ast}~w=w~\dot{\ast}~1=w$, for any $w\in\mathcal{A}^1_{\Lambda}$, 
\item[(ii)] $z_{k,s}w~\dot{\ast}~z_{l,t}w^{\prime}=z_{k+l,st}(w\ast w^{\prime})$ \\
for any $k,l\ge 1$ and any words $w,w^{\prime}\in\mathcal{A}^1_{\Lambda}$, 
\end{itemize}
and
\begin{itemize}
\item[(i)] $1~\dot{\overline{\ast}}~w=w~\dot{\overline{\ast}}~1=w$, for any $w\in\mathcal{A}^1_{\Lambda}$, 
\item[(ii)] $z_{k,s}w~\dot{\overline{\ast}}~z_{l,t}w^{\prime}=z_{k+l,st}(w~\overline{\ast}~w^{\prime})$ \\
for any $k,l\ge 1$ and any words $w,w^{\prime}\in\mathcal{A}^1_{\Lambda}$, 
\end{itemize}
respectively. } 
\end{defn}

As in \cite{H}, all of the four variations of harmonic products are shown to be associative and commutative. Note that the evaluation map $\mathtt{S}^{\ast}_{\bullet}$ is a $\overline{\ast}$-homomorphism and the evaluation map $\mathtt{s}^{\ast}_{\bullet}$ a $\dot{\overline{\ast}}$-homomorphism, i.e.,
\begin{equation}
\mathtt{S}^{\ast}_w(z)\mathtt{S}^{\ast}_{w^{\prime}}(z)=\mathtt{S}^{\ast}_{w\overline{\ast}w^{\prime}}(z),~\mathtt{s}^{\ast}_w(z)\mathtt{s}^{\ast}_{w^{\prime}}(z)=\mathtt{s}^{\ast}_{w\dot{\overline{\ast}}w^{\prime}}(z). \label{eq8}
\end{equation}

For $s\in\mathbb{C}^{\times}$, let $F_s$ denote $M_{\frac{1}{1-\delta(s)s}}\mathcal{I}^{-1}\iota\mathcal{I}M_{(1-\delta(s))+\delta(s)s}$. The following properties can be shown algebraically by induction. The proofs are a somewhat long and are given in $\S \ref{sec5}$. 
\begin{lem}\label{lem3-7}
{\it Let $\Lambda$ be a group. For any $w,w^{\prime}\in\mathcal{A}_{\Lambda}^1$, we have }

(i) {\it $d_{\ast}(w~\overline{\ast}~w^{\prime})=d_{\ast}(w)\ast d_{\ast}(w^{\prime})$, }

(i)${}^{\prime}$ {\it $d_{\ast}^{-1}(w\ast w^{\prime})=d_{\ast}^{-1}(w)~\overline{\ast}~d_{\ast}^{-1}(w^{\prime})$, }

(ii) {\it $d_{\ast}(w~\dot{\overline{\ast}}~w^{\prime})=d_{\ast}(w)~\dot{\ast}~d_{\ast}(w^{\prime})$, }

(ii)${}^{\prime}$ {\it $d_{\ast}^{-1}(w~\dot{\ast}~w^{\prime})=d_{\ast}^{-1}(w)~\dot{\overline{\ast}}~d_{\ast}^{-1}(w^{\prime})$. }
\end{lem}

\begin{lem}\label{lem3-8}
{\it Let $\Lambda$ be a subgroup of $\mathbb{C}^{\times}$ which contains an element $1-s$ for any $s \in \Lambda\backslash\{1\}$. For any $w\in\mathcal{A}_{\Lambda}^1$, any $w^{\prime}\in\mathcal{A}_{\{1\}}^1$, and any $s\in\mathbb{C}^{\times}$, we have }
$$ F_sd_{\ast}^{-1}(w)~\overline{\ast}~d_{\ast}^{-1}(w^{\prime})=F_sd_{\ast}^{-1}(w\ast w^{\prime}). $$
\end{lem}

\subsection{A Functional Equation}\label{subsec3-4}

\noindent We define the Newton series for the truncated MLV $\mathtt{s}_w^{\ast}(n)$ by
$$ \displaystyle \mathtt{s}_w^{\ast}(z)=\sum_{n=0}^{\infty}(-1)^n(\nabla \mathtt{s}_w^{\ast})(n)\binom{z}{n}. $$
If the subscript $w$ is a linear combination of words, the Newton series is also regarded as the corresponding linear combination of Newton series for each appearing word. We prove the following properties of convergence. (We have $\mathtt{s}_w^{\ast}(n)=\mathtt{s}_{\mathcal{I}(w)}^{\mathcyb{sh}}(n)$, and hence $(\nabla \mathtt{s}_w^{\ast})(n)=\mathtt{s}_{\star\mathcal{I}(w)}^{\mathcyb{sh}}(n)$ because of Theorem \ref{thm2-6}.)
\begin{prop}\label{prop3-3}
{\it Let $|s_i|\le 1$ for any $1\le i\le p$. Then, the Newton series
\begin{equation}
\displaystyle \sum_{n=0}^{\infty}(-1)^n\mathtt{s}^{\mathcyb{sh}}_{y_{s_1}\cdots y_{s_p}}(n)\binom{z}{n} \label{eq0-1}
\end{equation}
converges absolutely for $\mathrm{Re}(z)>-1$. If $s_1\neq 1$, (\ref{eq0-1}) converges for $\mathrm{Re}(z)>-2$.}
\end{prop}

\begin{proof}
Since $|\mathtt{s}^{\mathcyb{sh}}_{y_{s_1}\cdots y_{s_p}}(m)|=O((\log m)^{p-1}/m)$ as $m\to\infty$, we have the first assertion because of Corollary \ref{cor0-3-3}. Let $s_1\neq 1$. By Proposition \ref{prop3-1}, it is sufficient to show that the Dirichlet series
$$ \displaystyle \sum_{m=0}^{\infty}\frac{\mathtt{s}^{\mathcyb{sh}}_{y_{s_1}\cdots y_{s_p}}(m)}{(m+1)^{z+1}} $$
converges for $\mathrm{Re}(z)>-2$. Set $T(M)=\sum_{m=0}^M s_1^m=\frac{1-s_1^{M+1}}{1-s_1}$. Using the Abel summation method, we have
$$ \displaystyle \sum_{m=0}^M\frac{\mathtt{s}^{\mathcyb{sh}}_{y_{s_1}\cdots y_{s_p}}(m)}{(m+1)^{z+1}}=\sum_{M=m_1\ge\cdots m_p\ge 0}\frac{T(m_1-m_2)s_2^{m_2-m_3}\cdots s_{p-1}^{m_{p-1}-m_p}s_p^{m_p+1}}{(m_1+1)^{z+2}(m_2+1)\cdots (m_p+1)} $$
$$ \displaystyle +\sum_{M-1\ge m_1\ge\cdots\ge m_p\ge 0}\frac{T(m_1-m_2)s_2^{m_2-m_3}\cdots s_{p-1}^{m_{p-1}-m_p}s_p^{m_p+1}}{(m_2+1)\cdots (m_p+1)}\biggl\{\frac{1}{(m_1+1)^{z+2}}-\frac{1}{(m_1+2)^{z+2}}\biggl\}. $$
Since
$$ \displaystyle \sum_{M\ge m_2\ge\cdots\ge m_p\ge 0}\frac{1}{(m_2+1)\cdots (m_p+1)}=O((\log M)^{p-1})~~(m\to\infty), $$
there exists a constant $C_1$ such that
$$ \displaystyle \Biggl|\sum_{M=m_1\ge\cdots m_p\ge 0}\frac{T(m_1-m_2)s_2^{m_2-m_3}\cdots s_{p-1}^{m_{p-1}-m_p}s_p^{m_p+1}}{(m_1+1)^{z+2}(m_2+1)\cdots (m_p+1)}\Biggr| $$
$$ \displaystyle \le\frac{C_1}{(M+1)^{\sigma +2}}\sum_{M\ge m_2\ge\cdots\ge m_p\ge 0}\frac{1}{(m_2+1)\cdots (m_p+1)}\to 0~~(M\to\infty). $$
On the other hand, there exists a constant $C_2$ such that
$$ \displaystyle \sum_{m_1=0}^{\infty}\Biggl|\sum_{m_1\ge\cdots\ge m_p\ge 0}\frac{T(m_1-m_2)s_2^{m_2-m_3}\cdots s_{p-1}^{m_{p-1}-m_p}s_p^{m_p+1}}{(m_2+1)\cdots (m_p+1)}\biggl\{\frac{1}{(m_1+1)^{z+2}}-\frac{1}{(m_1+2)^{z+2}}\biggl\}\Biggr| $$
$$ \displaystyle \le C_2\sum_{m_1=0}^{\infty}\sum_{m_1\ge\cdots\ge m_p\ge 0}\frac{1}{(m_2+1)\cdots (m_p+1)}\Biggl|\frac{1}{(m_1+1)^{z+2}}-\frac{1}{(m_1+2)^{z+2}}\Biggr|. $$
Set $\sigma=\mathrm{Re}(z)$. We find that
$$ \displaystyle \Biggl|\frac{1}{(m_1+1)^{z+2}}-\frac{1}{(m_1+2)^{z+2}}\Biggr|=O\biggl(\frac{1}{m_1^{\sigma +2}}\biggr)~~(m_1\to\infty). $$
Hence,
$$ \displaystyle \sum_{m_1\ge\cdots\ge m_p\ge 0}\frac{1}{(m_2+1)\cdots (m_p+1)}\Biggl|\frac{1}{(m_1+1)^{z+2}}-\frac{1}{(m_1+2)^{z+2}}\Biggr|=O\biggl(\frac{(\log m_1)^{p-1}}{m_1^{\sigma +3}}\biggr) $$
as $m\to\infty$. Therefore, we conclude the Proposition.
\end{proof}

\begin{prop}\label{prop3-4}
{\it {\rm i)} Let $w=x^{k_1-1}y_{s_1}\cdots x^{k_n-1}y_{s_n},s_i(1\le i\le n)\in\mathbb{C}^{\times}$ with $|1-s_1\cdots s_i|\le 1$ for any $1\le i\le n$. Then, the Newton series $\mathtt{s}^{\ast}_w(z)$ converges absolutely for $\mathrm{Re}(z)>-1$. If $k_1=1$, $\mathtt{s}^{\ast}_w(z)$ converges for $\mathrm{Re}(z)>-2$. \\
{\rm ii)} Let $w=x^{k_1-1}y_{s_1}\cdots x^{k_n-1}y_{s_n},s_i(1\le i\le n)\in\mathbb{C}^{\times}$ with $|1-s_1\cdots s_i|\le 1$ for any $1\le i\le n$. Then, the Newton series $\mathtt{S}^{\ast}_w(z)$ converges absolutely for $\mathrm{Re}(z)>-1$. If $k_1=1$, $\mathtt{S}^{\ast}_w(z)$ converges for $\mathrm{Re}(z)>-2$. }
\end{prop}

\begin{proof}
The proof follows from Remark \ref{rmk2-7} 2 and Proposition \ref{prop3-3}. 
\end{proof}

\begin{prop}\label{prop0-3-7}
{\it Let $\Lambda=(0,1]\subset\mathbb{R}$, and $w=x^{k_1-1}y_{s_1}\cdots x^{k_n-1}y_{s_n}$ and $w^{\prime}=x^{k^{\prime}_1-1}y_{s^{\prime}_1}\cdots x^{k^{\prime}_{n^{\prime}}-1}y_{s^{\prime}_{n^{\prime}}}$ be words in $\mathcal{A}^1_{\Lambda}$. Then, the Newton series $\mathtt{s}_w^{\ast}(z)\mathtt{S}_{w^{\prime}}^{\ast}(z)$ is expressed as a Newton series that converges for $\mathrm{Re}(z)>-1$.}
\end{prop}

\begin{proof}
The proof follows from Remark \ref{rmk2-7} 2, (\ref{EQ}), Corollary \ref{cor4}, and Proposition \ref{prop0-3-6}. 
\end{proof}

\begin{thm}\label{thm3-9}
{\it Let $\Lambda=(0,1]\subset\mathbb{R}$, and $w=x^{k_1-1}y_{s_1}\cdots x^{k_n-1}y_{s_n}$ and $w^{\prime}=x^{k^{\prime}_1-1}y_{s^{\prime}_1}\cdots x^{k^{\prime}_{n^{\prime}}-1}y_{s^{\prime}_{n^{\prime}}}$ be words in $\mathcal{A}^1_{\Lambda}$, and $s\in\Lambda$. Then, we have
$$ \mathtt{s}^{\ast}_{L_{y_s}(w)}(z)\mathtt{S}^{\ast}_{w^{\prime}}(z)=\mathtt{s}^{\ast}_{L_{y_s}(w~\overline{\ast}~{w^{\prime}})}(z) $$
for $\mathrm{Re}(z)>-2$. }
\end{thm}

\begin{proof}
We find that each term of the right-hand side of the identity converges for $\mathrm{Re}(z)>-2$. By Proposition \ref{prop0-3-7}, the left-hand side is also expressed as a Newton series. Assuming Lemma \ref{lem3-8}, the identity holds on $\mathbb{Z}_{\ge 0}$. Therefore, the assertion is proven by Proposition \ref{prop3-2}.
\end{proof}

\subsection{Relations for MLV's}\label{subsec3-5}


For a positive integer $r$, we denote the set of $r$-th roots of unity by $\mu_r$. The MLV's $L^{\sharp}$ and $\overline{L}^{\sharp}$ ($\sharp=\mathcyr{sh}$ or $\ast$) were defined in $\S \ref{sec1}$.  We define the MLV-evaluation maps $\mathcal{L}^{\sharp}$, $\overline{\mathcal{L}}^{\sharp}:\mathcal{A}^0_{\mu_r}\to\mathbb{C}$ by $\mathbb{Q}$-linearity and
$$ \mathcal{L}^{\sharp}(x^{k_1-1}y_{s_1}\cdots x^{k_n-1}y_{s_n})=L^{\sharp}(k_1,\ldots ,k_n;r_1,\ldots ,r_n),\mathcal{L}^{\sharp}(1)=1 $$
and
$$ \overline{\mathcal{L}}^{\sharp}(x^{k_1-1}y_{s_1}\cdots x^{k_n-1}y_{s_n})=\overline{L}^{\sharp}(k_1,\ldots ,k_n;r_1,\ldots ,r_n),\overline{\mathcal{L}}^{\sharp}(1)=1, $$
respectively, where $s_i=\zeta^{r_i}$ ($1\le i\le n$) with $\zeta=\exp(2\pi i/r)$, which is a primitive $r$-th root of unity. We find that equation $\mathcal{L}^{\ast}(w)=\mathcal{L}^{\mathcyb{sh}}(\mathcal{I}(w))$ holds as well as equation (\ref{eq7}). 
\begin{prop}\label{prop3-10}
{\it Let $\Lambda$ be a subgroup of $\mathbb{C}^{\times}$ which contains an element $1-s$ for any $s \in \Lambda\backslash\{1\}$. In addition, let $w$ be a word of $\mathcal{A}^1_{\Lambda}$, and let $s$ be an element of $\Lambda$. For any non-negative integer $j$, we have
$$ \displaystyle \frac{1}{j!}\biggl(\frac{d}{dz}\biggr)^j \mathtt{s}^{\ast}_{\mathcal{I}^{-1}\iota L_{y_s}(w)}(-1)=(-1)^{j+1}\sum_{m=0}^{\infty}\mathtt{s}^{\ast}_{L_x^{-1}d_{\ast}^{-1}(\mathcal{I}^{-1}L_{x+\delta(s)y_s}\varphi d_{\mathcyb{sh}}(w)~\dot{\ast}~y_1^{j+1})}(m). $$
In particular, if $\Lambda=\mu_r$ for a positive integer $r$, we have
$$ \displaystyle \frac{1}{j!}\biggl(\frac{d}{dz}\biggr)^j \mathtt{s}^{\ast}_{\mathcal{I}^{-1}\iota L_{y_s}d_{\mathcyb{sh}}^{-1}\mathcal{I}M_s(w)}(-1)=(-1)^{j+1}\mathcal{L}^{\mathcyb{sh}}(L_x^{-1}\mathcal{I}(\mathcal{I}^{-1}L_{x+\delta(s)y_s}\varphi \mathcal{I}M_s(w)~\dot{\ast}~y_1^{j+1}). $$ }
\end{prop}

\begin{proof}
Generally, $j$-th coefficients of Taylor expansion at $z=-1$ of Newton series $a(z)=\sum_{n=0}^{\infty}(-1)^j(\nabla a)(n)\binom{z}{n}$ can be expressed as
$$ \displaystyle \frac{1}{j!}a^{(j)}(-1)=(-1)^j\sum_{n=0}^{\infty}(\nabla a)(n)\mathtt{S}^{\ast}_{d_{\ast}^{-1}(y_1^j)}(n-1) $$
(see \cite[Lemma 4.6]{K}). Using the identities 
$$ \mathtt{s}^{\ast}_{\mathcal{I}^{-1}(w)}=\mathtt{s}^{\mathcyb{sh}}_w(z) $$
and
$$ \mathtt{S}^{\ast}_{d_{\ast}^{-1}(y_1^j)}(n-1)=(n+1)\mathtt{s}^{\ast}_{d_{\ast}^{-1}(y_1^{j+1})}(n) $$
and Theorem \ref{thm2-6}, we have
\begin{eqnarray}
&{}& \displaystyle \frac{1}{j!}\biggl(\frac{d}{dz}\biggr)^j \mathtt{s}^{\ast}_{\mathcal{I}^{-1}\iota L_{y_s}(w)}(-1)=\frac{1}{j!}\biggl(\frac{d}{dz}\biggr)^j \mathtt{s}^{\mathcyb{sh}}_{\iota L_{y_s}(w)}(-1) \nonumber \\
&=&(-1)^j\sum_{n=0}^{\infty}\mathtt{s}^{\mathcyb{sh}}_{\star\iota L_{y_s}(w)}(n)\cdot (n+1)\mathtt{s}^{\ast}_{d_{\ast}^{-1}(y_1^{j+1})}(n). \label{eq9}
\end{eqnarray}
Note that $\star\iota=-d_{\mathcyb{sh}}^{-1}\varphi d_{\mathcyb{sh}}$ because of Lemma \ref{lem2-3} and $d_{\mathcyb{sh}}\iota=\iota d_{\mathcyb{sh}}$. Using
$$ \mathtt{s}^{\mathcyb{sh}}_{\star\iota L_{y_s}(w)}(n)=\mathtt{s}^{\ast}_{\mathcal{I}^{-1}\star\iota L_{y_s}(w)}(n)=-\mathtt{s}^{\ast}_{\mathcal{I}^{-1}d_{\mathcyb{sh}}^{-1}\varphi d_{\mathcyb{sh}}L_{y_s}(w)}(n) $$
and equation (\ref{eq8}), we obtain
\begin{equation}
(\ref{eq9})=(-1)^{j+1}\sum_{n=0}^{\infty}(n+1)\mathtt{s}^{\ast}_{\mathcal{I}^{-1}d_{\mathcyb{sh}}^{-1}\varphi d_{\mathcyb{sh}}L_{y_s}(w)~\dot{\overline{\ast}}~d_{\mathcyb{sh}}^{-1}(y_1^{j+1})}(n). \label{eq10}
\end{equation}
Since $\mathcal{I}^{-1}d_{\mathcyb{sh}}^{-1}\varphi d_{\mathcyb{sh}}L_{y_s}(w)~\dot{\overline{\ast}}~d_{\mathcyb{sh}}^{-1}(y_1^{j+1})\in\mathcal{A}^0_{\Lambda}$, we have
$$ (n+1)\mathtt{s}^{\ast}_{\mathcal{I}^{-1}d_{\mathcyb{sh}}^{-1}\varphi d_{\mathcyb{sh}}L_{y_s}(w)~\dot{\overline{\ast}}~d_{\mathcyb{sh}}^{-1}(y_1^{j+1})}(n)=\mathtt{s}^{\ast}_{L_x^{-1}(\mathcal{I}^{-1}d_{\mathcyb{sh}}^{-1}\varphi d_{\mathcyb{sh}}L_{y_s}(w)~\dot{\overline{\ast}}~d_{\mathcyb{sh}}^{-1}(y_1^{j+1}))}(n). $$
By this identity, $\mathcal{I}^{-1}d_{\mathcyb{sh}}^{-1}=d_{\ast}^{-1}\mathcal{I}^{-1}$, Lemma \ref{lem3-7} (i)${}^{\prime}$, and $\varphi d_{\mathcyb{sh}}L_{y_s}=L_{x+\delta(s)y_s}\varphi d_{\mathcyb{sh}}$, we have the first assertion
$$ (\ref{eq10})=(-1)^{j+1}\sum_{m=0}^{\infty}\mathtt{s}^{\ast}_{L_x^{-1}d_{\ast}^{-1}(\mathcal{I}^{-1}L_{x+\delta(s)y_s}\varphi d_{\mathcyb{sh}}(w)~\dot{\ast}~y_1^{j+1})}(m). $$
Substituting $d_{\mathcyb{sh}}^{-1}\mathcal{I}M_s(w)$ into $w$ in the above equality, we have
\begin{equation}
\displaystyle \frac{1}{j!}\biggl(\frac{d}{dz}\biggr)^j \mathtt{s}^{\ast}_{\mathcal{I}^{-1}\iota L_{y_s}d_{\mathcyb{sh}}^{-1}\mathcal{I}M_s(w)}(-1)=(-1)^{j+1}\sum_{m=0}^{\infty}\mathtt{s}^{\ast}_{L_x^{-1}d_{\ast}^{-1}(\mathcal{I}^{-1}L_{x+\delta(s)y_s}\varphi\mathcal{I}M_s(w)~\dot{\ast}~y_1^{j+1})}(m). \label{eq11}
\end{equation}
If $\Lambda=\mu_r$, every term of $L_x^{-1}d_{\ast}^{-1}(\mathcal{I}^{-1}L_{x+\delta(s)y_s}\varphi\mathcal{I}M_s(w)~\dot{\ast}~y_1^{j+1})$ is admissible and converges to a linear combination of MLV's, i.e., 
\begin{eqnarray*}
(\ref{eq11})&=&(-1)^{j+1}\overline{\mathcal{L}}^{\ast}(L_x^{-1}d_{\ast}^{-1}(\mathcal{I}^{-1}L_{x+\delta(s)y_s}\varphi\mathcal{I}M_s(w)~\dot{\ast}~y_1^{j+1})) \\
&=&(-1)^{j+1}\mathcal{L}^{\mathcyb{sh}}(L_x^{-1}\mathcal{I}(\mathcal{I}^{-1}L_{x+\delta(s)y_s}\varphi\mathcal{I}M_s(w)~\dot{\ast}~y_1^{j+1})). 
\end{eqnarray*}
\end{proof}

\begin{rmk}\label{rmk3-11}
Let $w\in\mathcal{A}^1_{\{1\}}$. Differentiating both sides of the identity
$$ \mathtt{S}^{\ast}_w(z)=(z+1)\mathtt{s}^{\ast}_{L_{y_1}(w)}(z), $$
we see that
$$ \displaystyle \frac{1}{j!}\biggl(\frac{d}{dz}\biggr)^j\mathtt{S}^{\ast}_{d_{\ast}^{-1}(w)}(-1)=\left\{
\begin{array}{ll}
0 & j=0, \\
(-1)^j \mathcal{L}^{\mathcyb{sh}}(L_x^{-1}(L_x\varphi(w^{\prime})~\dot{\ast}~y_1^j)) & j>0.
\end{array}
\right. $$
\end{rmk}

Let $\mathcal{A}^1_{\Lambda,>0}$ denote the set of polynomials in $\mathcal{A}_{\Lambda}^1$ without constant terms. 
\begin{thm}\label{thm3-12}
{\it For any $w\in\mathcal{A}^1_{\mu_r,>0}$, any $w^{\prime}\in\mathcal{A}^1_{\{1\},>0}$, and any $m\ge 0$, we have
\begin{flushleft}
$ \displaystyle \sum_{k+l=m, k\ge 0,l>0}\mathcal{L}^{\mathcyb{sh}}(L_x^{-1}\mathcal{I}(\mathcal{I}^{-1}L_{x+\delta(s)y_s}\varphi \mathcal{I}M_s(w)~\dot{\ast}~y_1^{k+1}))\mathcal{L}^{\mathcyb{sh}}(L_x^{-1}(L_x\varphi(w^{\prime})~\dot{\ast}~y_1^l)) \qquad $
\end{flushleft}
\begin{flushright}
$ \displaystyle =\mathcal{L}^{\mathcyb{sh}}(L_x^{-1}\mathcal{I}(\mathcal{I}^{-1}L_{x+\delta(s)y_s}\varphi \mathcal{I}M_s(w\ast w^{\prime})~\dot{\ast}~y_1^{m+1})). $
\end{flushright}
If $m=0$, we consider the left-hand side as $0$. }
\end{thm}

\begin{proof}
First, we find that
\begin{equation}
\left.
\begin{array}{ll}
\displaystyle \sum_{\begin{subarray}{c}k+l=m, \\ k\ge 0,l>0\end{subarray}}\mathrm{Li}^{\mathcyb{sh}}_{L_x^{-1}\mathcal{I}(\mathcal{I}^{-1}L_{x+\delta(s)y_s}\varphi \mathcal{I}M_s(w)~\dot{\ast}~y_1^{k+1})}(1)\mathrm{Li}^{\mathcyb{sh}}_{L_x^{-1}(L_x\varphi(w^{\prime})~\dot{\ast}~y_1^l)}(1)\qquad \\
\qquad\qquad \displaystyle =\mathrm{Li}^{\mathcyb{sh}}_{L_x^{-1}\mathcal{I}(\mathcal{I}^{-1}L_{x+\delta(s)y_s}\varphi \mathcal{I}M_s(w\ast w^{\prime})~\dot{\ast}~y_1^{m+1})}(1)
\end{array}
\right. \label{id5}
\end{equation}
for any $w\in\mathcal{A}^1_{(0,1],>0}$, any $w^{\prime}\in\mathcal{A}^1_{\{1\},>0}$, and any $m\ge 0$. Since each term in (\ref{id5}) is holomorphic for $|s_i|<1$ ($1\le i\le n$) because of the Identity Theorem. Moreover, we see that each term in (\ref{id5}) is continuous for $|s_i|\le 1$ ($1\le i\le n$), hence the identity (\ref{id5}) is valid for $|s_i|\le 1$ ($1\le i\le n$). If $s_i\in\mu_r$, each value of non-strict MPL's at $z=1$ appearing in (\ref{id5}) is nothing but the MLV. Replacing $w$ with $d_{\mathcyb{sh}^{-1}}\mathcal{I}M_s(w)$, Theorem is proven. 
\end{proof}

Finally, we give the case of $m=0$, the linear part of Theorem \ref{thm3-12}, which is shown to contain the extended derivation relation in the next section. 
\begin{cor}\label{cor3-13}
{\it Let $r$ be a positive integer. For any $s\in\mu_r$, we have 
$$ L_{x+\delta(s)y_s}\varphi \mathcal{I}M_s(\mathcal{A}^1_{\mu_r,>0}\ast\mathcal{A}^1_{\{1\},>0})\subset\mathrm{Ker}\mathcal{L}^{\mathcyb{sh}}.$$}
\end{cor}


\section{Extended Derivation}\label{sec4}

\noindent A generalization of the derivation relation was first mentioned in \cite{IKZ} for the multiple zeta value (MZV) case, formulated as a conjecture in \cite{Kan}, and proven in \cite{T} by reducing the relation to the relations found in \cite{K}. In this section, we present a generalization of the derivation relation for the MLV case, which has been suggested to the authors by Masanobu Kaneko. It is worthwhile to verify several advantageous properties of the extended derivation operators $\widehat{\partial}_n^{(c)}$ and $\partial_n^{(c)}$. 

Throughout this section, we fix a positive integer $r$. As defined in the previous section, let $\mu_r=\{1,\zeta,\zeta^2,\ldots,\zeta^{r-1}\}$ with $\zeta=\exp(2\pi i/r)$. In addition, we denote $x+y_1(\in\mathcal{A}_{\mu_r})$ by $z$ in this section. 

\subsection{Definitions and Properties}\label{subsec4-1}

\noindent For $n\ge 1$, the derivation operator $\partial_n:\mathcal{A}_{\mu_r}\to\mathcal{A}_{\mu_r}$ appeared in \cite{AK} satisfies
$$ \displaystyle \partial_n=\frac{1}{(n-1)!}\mathrm{ad}(\theta)^{n-1}(\partial_1), $$
where $\theta$ is a derivation on $\mathcal{A}_{\mu_r}$ defined by $\theta(u)=\frac{1}{2}(uz+zu)$ for $u=x$ or $y_s$ ($s\in\mu_r$), $\partial_1:\mathcal{A}_{\mu_r}\to\mathcal{A}_{\mu_r}$ a derivation characterized by 
$$ \partial_1(x)=xy_1,~\partial_1(y_s)=-xy_s+y_sy_1-y_sy_s, $$
and $\mathrm{ad}(\theta)(\partial)=[\theta, \partial]:=\theta\partial -\partial\theta$. Then, we give an extension of $\partial_n$ as follows. 
\begin{defn}\label{defn4-1}
{\it Let $c$ be a rational number, and let $H$ be the derivation on $\mathcal{A}_{\mu_r}$ defined by $H(w)=\deg (w)w$ for any words $w \in \mathcal{A}_{\mu_r}$. For each integer $n \ge 1$, we define a $\mathbb{Q}$-linear map $\widehat{\partial}_n^{(c)}$ from $\mathcal{A}_{\mu_r}$ to $\mathcal{A}_{\mu_r}$ by}
$$ \displaystyle \widehat{\partial}_n^{(c)}=\frac{1}{(n-1)!}\mathrm{ad}(\widehat{\theta}^{(c)})^{n-1}(\partial_1), $$
{\it where $\widehat{\theta}^{(c)}$ is the $\mathbb{Q}$-linear map defined by $\widehat{\theta}^{(c)}(u)=\theta (u)$ ($u=x$ or $y_s$) and the rule}
$$ \widehat{\theta}^{(c)}(ww^{\prime})=\widehat{\theta}^{(c)}(w)w^{\prime}+w\widehat{\theta}^{(c)}(w^{\prime})+cH(w)\partial_1(w^{\prime}) $$
{\it for any $w,w^{\prime} \in \mathcal{A}_{\mu_r}$.}
\end{defn}

If $c=0$, the extended derivation $\widehat{\partial}_n^{(c)}$ is reduced to the ordinary derivation $\partial_n$. It is known in \cite{AK} that $\partial_n(\mathcal{A}_{\mu_r}^0)\subset\mathrm{Ker}\mathcal{L}^{\mathcyb{sh}}$ holds for any $n\ge 1$, which is called the derivation relation. Although, if $c\neq 0$ and $n\ge 2$, the operator $\widehat{\partial}_n^{(c)}$ is no longer a derivation, and we cannot prove (but can expect experimentally) that the extended derivation relation is contained in the EDSR, we find in Theorem \ref{thm4-14} or (\ref{eq}) that the extended derivation relation is proven by using Corollary \ref{cor3-13}. 

Now, we define the operator $\widehat{\psi}_n^{(c)}(u)$, which is closely related to the operator $\widehat{\partial}_n^{(c)}$, and then show some properties of this operator. 
\begin{defn}\label{defn4-2}
{\it Let $u\in\mathbb{Q}\cdot x+\sum_{s\in\mu_r}\mathbb{Q}\cdot y_s$. For $n\ge 1$ and $c\in\mathbb{Q}$, we define the sequence of operators $\{\widehat{\psi}_n^{(c)}(u)\}_{n=1}^{\infty}$ by $\widehat{\psi}_1^{(c)}(u)=L_{\partial_1(u)}$ and
$$ \displaystyle \widehat{\psi}_n^{(c)}(u)=\frac{1}{n-1}\bigl([\widehat{\theta}^{(c)},\widehat{\psi}_{n-1}^{(c)}(u)]-\frac{1}{2}(L_z\widehat{\psi}_{n-1}^{(c)}(u)+\widehat{\psi}_{n-1}^{(c)}(u)L_z)-c\widehat{\psi}_{n-1}^{(c)}(u)\partial_1\bigl) $$
for $n\ge 2$ (where $L_w$ has been defined in $\S \ref{subsec3-3}$ by $L_w(w^{\prime})=ww^{\prime}$). }
\end{defn}
\noindent Let $\nu$ be a $\mathbb{Q}$-linear operator on $\mathbb{Q}\cdot x+\sum_{s\in\mu_r}\mathbb{Q}\cdot y_s$ given by $\nu(x)=0$ and $\nu(y_s)=1$ for any $s\in\mu_r$.
\begin{lem}\label{lem4-3}
{\it Let $n\ge 1$ and $c\in\mathbb{Q}$. There is another sequence of operators $\{\widehat{\phi}_n^{(c)}\}_{n=0}^{\infty}$ such that
$$ \widehat{\psi}_n^{(c)}(u)=(-1)^{\nu(u)}\bigl(L_x\widehat{\phi}_{n-1}^{(c)}L_{y_1+\nu(u)(u-y_1)}+\nu(u)L_u\widehat{\phi}_{n-1}^{(c)}L_{u-y_1}\bigr) $$
where $u\in\mathbb{Q}\cdot x+\sum_{s\in\mu_r}\mathbb{Q}\cdot y_s$. }
\end{lem}
\noindent The proof of Lemma \ref{lem4-3} is given in $\S.\ref{sec5}$. 
\begin{cor}\label{cor4-4}
{\it For $n\ge 1$ and $c\in\mathbb{Q}$, we have }
\begin{itemize}
\item[(i)]$\widehat{\psi}_n^{(c)}(x+y_1)=0.$
\item[(ii)]$\widehat{\psi}_n^{(c)}(x+\delta(s)y_s)=L_{x+\delta(s)y_s}\widehat{\phi}_{n-1}^{(c)}L_{y_1-\delta(s)y_s}~(s\in\mu_r).$
\end{itemize}
\end{cor}
\begin{proof}
By Lemma \ref{lem4-3}, we see that
$$ \widehat{\psi}_n^{(c)}(x)=L_x\widehat{\phi}_n^{(c)}L_{y_1},~\widehat{\psi}_n^{(c)}(y_1)=-L_x\widehat{\phi}_n^{(c)}L_{y_1}, $$
which implies (i) and (ii) for $s=1$. If $s\neq 1$, then
$$ \widehat{\psi}_n^{(c)}(y_s)=-L_x\widehat{\phi}_n^{(c)}L_{y_s}+L_{y_s}\widehat{\phi}_n^{(c)}L_{y_s-y_1} $$
by Lemma \ref{lem4-3} and we obtain (ii) for $s\neq 1$. 
\end{proof}

Next, we present the commutativity of $\widehat{\partial}_n^{(c)}$. 
\begin{prop}\label{prop4-5}
{\it For any $n,m\ge 1$ and any $c,c^{\prime}\in\mathbb{Q}$, we have $[\widehat{\partial}_n^{(c)},\widehat{\partial}_m^{(c^{\prime})}]=0$. }
\end{prop}
\begin{proof}
Let $n\ge 1$ and $u\in\mathbb{Q}\cdot x+\sum_{s\in\mu_r}\mathbb{Q}\cdot y_s$. We prove the following statements $(\mathrm{A}_n)$ and $(\mathrm{B}_n)$ inductively as $(\mathrm{A}_1),(\mathrm{B}_1)\Rightarrow(\mathrm{A}_2)\Rightarrow(\mathrm{B}_2)\Rightarrow(\mathrm{A}_3)\Rightarrow\cdots$. 
\begin{itemize}
\item[$(\mathrm{A}_n)$]$[\widehat{\partial}_n^{(c)},L_u]=\widehat{\psi}_n^{(c)}(u).$
\item[$(\mathrm{B}_n)$]$[\widehat{\partial}_n^{(c)},\widehat{\partial}_i^{(c^{\prime})}]=0$, for any $1\le i\le n$, and any $c,c^{\prime}\in\mathbb{Q}.$
\end{itemize}
If $(\mathrm{B}_n)$ is proven, we obtain the assertion. 

Here, we present three remarks. First, if $(\mathrm{A}_n)$ is proven, we see
\begin{itemize}
\item[$(\alpha_n)$]$[\widehat{\partial}_n^{(c)},L_z]=0$
\end{itemize}
because of Corollary \ref{cor4-4}(i). 

Second, let 
\begin{itemize}
\item[$(\mathrm{B}_{n,i})$]$[\widehat{\partial}_n^{(c)},\widehat{\partial}_i^{(c^{\prime})}]=0$ for a fixed $1\le i\le n$ and any $c,c^{\prime}\in\mathbb{Q}$
\end{itemize}
for $1\le i\le n$. Showing statement $(\mathrm{B}_n)$ is equivalent to showing statement $(\mathrm{B}_{n,i})$ for all $1\le i\le n$. Moreover, setting
\begin{itemize}
\item[$(\mathrm{B}_{n,i}^{\prime})$]$[[\widehat{\partial}_n^{(c)},\widehat{\partial}_i^{(c^{\prime})}],L_u]=0$ for a fixed $1\le i\le n$ and any $c,c^{\prime}\in\mathbb{Q}$,
\end{itemize}
we see that each $(\mathrm{B}_{n,i})$ is equivalent to $(\mathrm{B}_{n,i}^{\prime})$. Instead of $(\mathrm{B}_n)$, we show $(\mathrm{B}_{n,i}^{\prime})$ by induction on $i$. 

Third, when the statement $(\mathrm{A}_i)$ (hence $(\alpha_i)$) and the statement $(\mathrm{B}_i)$ are proven for any $1\le i\le n$, we can consider the commutative polynomial ring $\mathbb{Q}[L_z,\widehat{\partial}_1^{(c)},\ldots,\widehat{\partial}_n^{(c)}]$. Define the degree of generators by $\deg(L_z)=1$ and $\deg(\widehat{\partial}_d^{(c)})=d$. Let $\mathbb{Q}[L_z,\widehat{\partial}_1^{(c)},\ldots,\widehat{\partial}_n^{(c)}]_{(d)}$ denote the degree-$d$ homogenous part. Then, we obtain, 
\begin{itemize}
\item[$(\beta_n)$]$\widehat{\phi}_n^{(c)}\in\mathbb{Q}[L_z,\widehat{\partial}_1^{(c)},\ldots,\widehat{\partial}_n^{(c)}]_{(n)}$
\end{itemize}
because of the recursive rule (\ref{eq14}) of the operator $\widehat{\phi}_n^{(c)}$. 

We now begin to prove $(\mathrm{A}_n)$ and $(\mathrm{B}_n)$. Since $\widehat{\partial}_1^{(c)}=\partial_1$ for any $c\in\mathbb{Q}$ and 
$$ [\partial_1,L_u](w)=\partial_1(uw)-u\partial_1(w)=\partial_1(u)w=L_{\partial_1(u)}(w) $$
for $w\in\mathcal{A}_{\mu_r}$, the statement $(\mathrm{A}_1)$ holds. The statement $(\mathrm{B}_1)$ is trivial because $\widehat{\partial}_1^{(c)}=\partial_1^{(c)}$ for any $c\in\mathbb{Q}$. 

Suppose $(\mathrm{A}_n)$ (hence also $(\alpha_n)$) and $(\mathrm{B}_n)$ are proven. Then, 
\begin{eqnarray*}
\displaystyle n[\widehat{\partial}_{n+1}^{(c)},L_u] &=& [[\widehat{\theta}^{(c)},\widehat{\partial}_n^{(c)}],L_u] \\
\displaystyle &=& -[[\widehat{\partial}_n^{(c)}],L_u],\widehat{\theta}^{(c)}]-[[L_u,\widehat{\theta}^{(c)}],\widehat{\partial}_n^{(c)}]] \\
\displaystyle &=& [\widehat{\theta}^{(c)},\widehat{\psi}_n^{(c)}(u)]+[L_{\widehat{\theta}^{(c)}(u)}+cL_u\partial_1,\widehat{\partial}_n^{(c)}] \\
\displaystyle &=& [\widehat{\theta}^{(c)},\widehat{\psi}_n^{(c)}(u)]-[\widehat{\partial}_n^{(c)},\frac{1}{2}(L_uL_z+L_zL_u)+cL_u\partial_1] \\
\displaystyle &=& [\widehat{\theta}^{(c)},\widehat{\psi}_n^{(c)}(u)]-\frac{1}{2}(L_z\widehat{\psi}_n^{(c)}(u)+\widehat{\psi}_n^{(c)}(u)L_z)-c\widehat{\psi}_n^{(c)}(u)\partial_1.
\end{eqnarray*}
Hence, we obtain $[\widehat{\partial}_{n+1}^{(c)},L_u]=\widehat{\psi}_{n+1}^{(c)}(u)$, and $(\mathrm{A}_{n+1})$ (hence $(\alpha_{n+1})$) holds. 

Next, suppose that all of $(\mathrm{A}_j)$'s (hence $(\alpha_j)$'s) $(1\le j\le n+1)$ and all of $(\mathrm{B}_j)$'s (hence $(\beta_j)$'s) $(1\le j\le n)$ are proven. We show $(\mathrm{B}_{n+1,i}^{\prime})$ $(1\le j\le n+1)$ by induction on $i$. 
\begin{eqnarray*}
\displaystyle &{}& [[\widehat{\partial}_{n+1}^{(c)},\widehat{\partial}_i^{(c^{\prime})}],L_u] \\
\displaystyle &=& -[[\widehat{\partial}_i^{(c^{\prime})},L_u],\widehat{\partial}_{n+1}^{(c)}]-[[L_u,\widehat{\partial}_{n+1}^{(c)}],\widehat{\partial}_i^{(c^{\prime})}] \\
\displaystyle &=& [\widehat{\partial}_{n+1}^{(c)},\widehat{\psi}_i^{(c^{\prime})}(u)]-[\widehat{\partial}_i^{(c^{\prime})},\widehat{\psi}_{n+1}^{(c)}(u)] \\
\displaystyle &=& (-1)^{\nu(u)}\bigl([\widehat{\partial}_{n+1}^{(c)},L_x\widehat{\phi}_{i-1}^{(c^{\prime})}L_{y_1+\nu(u)(u-y_1)}+\nu(u)L_u\widehat{\phi}_{i-1}^{(c^{\prime})}L_{u-y_1}] \\
\displaystyle &{}& \quad -[\widehat{\partial}_i^{(c^{\prime})},L_x\widehat{\phi}_n^{(c)}L_{y_1+\nu(u)(u-y_1)}+\nu(u)L_u\widehat{\phi}_n^{(c)}L_{u-y_1}]\bigl).
\end{eqnarray*}
We must show that this becomes $0$. Here, set

\quad $ P=[\widehat{\partial}_{n+1}^{(c)},L_x\widehat{\phi}_{i-1}^{(c^{\prime})}L_{y_1+\nu(u)(u-y_1)}]-[\widehat{\partial}_i^{(c^{\prime})},L_x\widehat{\phi}_n^{(c)}L_{y_1+\nu(u)(u-y_1)}], $

\quad $ Q=[\widehat{\partial}_{n+1}^{(c)},L_u\widehat{\phi}_{i-1}^{(c^{\prime})}L_{u-y_1}]-[\widehat{\partial}_i^{(c^{\prime})},L_u\widehat{\phi}_n^{(c)}L_{u-y_1}]. $ \\
Then, we need only show that 
\begin{equation}
\nu(u)Q=-P. \label{eq15}
\end{equation}
Note that $(-1)^{\nu(u)}\nu(u)=-\nu(u)$. 

Suppose $i=1$. Since $\widehat{\phi}_{i-1}^{(c^{\prime})}=\mathrm{id}_{\mathcal{A}_{\mu_r}},~\widehat{\partial}_i^{(c^{\prime})}=\partial_1$, we have
$$ [\widehat{\partial}_{n+1}^{(c)},\widehat{\phi}_{i-1}^{(c^{\prime})}]=0,~[\widehat{\partial}_i^{(c^{\prime})},\widehat{\phi}_n^{(c)}]=0. $$
Hence,
\begin{eqnarray*}
\displaystyle P &=& \widehat{\psi}_{n+1}^{(c)}(x)L_{y_1+\nu(u)(u-y_1)}+L_x\widehat{\psi}_{n+1}^{(c)}(y_1+\nu(u)(u-y_1)) \\
\displaystyle &{}& \quad -\widehat{\psi}_i^{(c^{\prime})}(x)\widehat{\phi}_n^{(c)}L_{y_1+\nu(u)(u-y_1)}-L_x\widehat{\phi}_n^{(c)}\widehat{\psi}_1^{(c^{\prime})}(y_1+\nu(u)(u-y_1)) \\
\displaystyle &=& L_x\widehat{\phi}_n^{(c)}L_{y_1}L_{y_1+\nu(u)(u-y_1)}-L_xL_x\widehat{\phi}_n^{(c)}L_{y_1}+\nu(u)L_xL_x\widehat{\phi}_n^{(c)}L_{y_1} \\
\displaystyle &{}& \quad +(-1)^{\nu(u)}\nu(u)L_x(L_x\widehat{\phi}_n^{(c)}L_{y_1+\nu(u)(u-y_1)}+\nu(u)L_u\widehat{\phi}_n^{(c)}L_{u-y_1}) \\
\displaystyle &{}& \quad -L_xL_{y_1}\widehat{\phi}_n^{(c)}L_{y_1+\nu(u)(u-y_1)}+L_x\widehat{\phi}_n^{(c)}L_xL_{y_1}-\nu(u)L_x\widehat{\phi}_n^{(c)}L_xL_{y_1} \\
\displaystyle &{}& \quad -(-1)^{\nu(u)}\nu(u)L_x\widehat{\phi}_n^{(c)}(L_xL_{y_1+\nu(u)(u-y_1)}+\nu(u)L_uL_{u-y_1}) \\
\displaystyle &=& L_x[\widehat{\phi}_n^{(c)},L_{y_1}]L_{y_1+\nu(u)(u-y_1)}+L_x[\widehat{\phi}_n^{(c)},L_x]L_{y_1}-\nu(u)L_x[\widehat{\phi}_n^{(c)},L_x]L_{y_1} \\
\displaystyle &{}& \quad +\nu(u)L_x[\widehat{\phi}_n^{(c)},L_x]L_{y_1+\nu(u)(u-y_1)}+\nu(u)L_x[\widehat{\phi}_n^{(c)},L_u]L_{u-y_1} \\
\displaystyle &=& \nu(u)L_x[\widehat{\phi}_n^{(c)},L_u]L_{u-y_1}.
\end{eqnarray*}
On the other hand, 
\begin{eqnarray*}
\displaystyle Q &=& \widehat{\psi}_{n+1}^{(c)}(u)L_{u-y_1}+L_u\widehat{\psi}_{n+1}^{(c)}(u-y_1)-\widehat{\psi}_1^{(c^{\prime})}(u)\widehat{\phi}_n^{(c)}L_{u-y_1}-L_u\widehat{\phi}_n^{(c)}\widehat{\psi}_1^{(c^{\prime})}(u-y_1) \\
\displaystyle &=& (-1)^{\nu(u)}\bigl(L_x\widehat{\phi}_n^{(c)}L_{y_1+\nu(u)(u-y_1)}+\nu(u)L_u\widehat{\phi}_n^{(c)}L_{u-y_1}\bigl)L_{u-y_1} \\
\displaystyle &{}& \quad +(-1)^{\nu(u)}L_u\bigl(L_x\widehat{\phi}_n^{(c)}L_{y_1+\nu(u)(u-y_1)}+\nu(u)L_u\widehat{\phi}_n^{(c)}L_{u-y_1}\bigl)+L_uL_x\widehat{\phi}_n^{(c)}L_{y_1} \\
\displaystyle &{}& \quad -(-1)^{\nu(u)}(L_xL_{y_1+\nu(u)(u-y_1)}+\nu(u)L_uL_{u-y_1})\widehat{\phi}_n^{(c)}L_{u-y_1} \\
\displaystyle &{}& \quad -(-1)^{\nu(u)}L_u\widehat{\phi}_n^{(c)}(L_xL_{y_1+\nu(u)(u-y_1)}+\nu(u)L_uL_{u-y_1})-L_u\widehat{\phi}_n^{(c)}L_xL_{y_1} \\
\displaystyle &=& (-1)^{\nu(u)}L_x[\widehat{\phi}_n^{(c)},L_{y_1+\nu(u)(u-y_1)}]L_{u-y_1}-\nu(u)L_u[\widehat{\phi}_n^{(c)},L_{u-y_1}]L_{u-y_1} \\
\displaystyle &{}& \quad -(-1)^{\nu(u)}L_u[\widehat{\phi}_n^{(c)},L_x]L_{y_1+\nu(u)(u-y_1)}+\nu(u)L_u[\widehat{\phi}_n^{(c)},L_u]L_{u-y_1}-L_u[\widehat{\phi}_n^{(c)},L_x]L_{y_1}.
\end{eqnarray*}
Thus, we obtain $\nu(u)Q=-\nu(u)L_x[\widehat{\phi}_n^{(c)},L_u]L_{u-y_1}=-P$ and (\ref{eq15}) is proven. 

Suppose $i>1$. By the induction hypothesis, again, 
$$[\widehat{\partial}_{n+1}^{(c)},\widehat{\phi}_{i-1}^{(c^{\prime})}]=0,~[\widehat{\partial}_i^{(c^{\prime})},\widehat{\phi}_n^{(c)}]=0, $$
and we have
\begin{eqnarray*}
\displaystyle P &=& \widehat{\psi}_{n+1}^{(c)}(x)\widehat{\phi}_{i-1}^{(c^{\prime})}L_{y_1+\nu(u)(u-y_1)}+L_x\widehat{\phi}_{i-1}^{(c^{\prime})}\widehat{\psi}_{n+1}^{(c)}(y_1+\nu(u)(u-y_1)) \\
\displaystyle &{}& \quad -\widehat{\psi}_i^{(c^{\prime})}(x)\widehat{\phi}_n^{(c)}L_{y_1+\nu(u)(u-y_1)}-L_x\widehat{\phi}_n^{(c)}\widehat{\psi}_1^{(c^{\prime})}(y_1+\nu(u)(u-y_1)) \\
\displaystyle &=& L_x\widehat{\phi}_n^{(c)}L_{y_1}\widehat{\phi}_{i-1}^{(c^{\prime})}L_{y_1+\nu(u)(u-y_1)}-L_x\widehat{\phi}_{i-1}^{(c^{\prime})}L_x\widehat{\phi}_n^{(c)}L_{y_1}+\nu(u)L_x\widehat{\phi}_{i-1}^{(c^{\prime})}L_x\widehat{\phi}_n^{(c)}L_{y_1} \\
\displaystyle &{}& \quad +(-1)^{\nu(u)}\nu(u)L_x\widehat{\phi}_{i-1}^{(c^{\prime})}(L_x\widehat{\phi}_n^{(c)}L_{y_1+\nu(u)(u-y_1)}+\nu(u)L_u\widehat{\phi}_n^{(c)}L_{u-y_1}) \\
\displaystyle &{}& \quad -L_x\widehat{\phi}_{i-1}^{(c^{\prime})}L_{y_1}\widehat{\phi}_n^{(c)}L_{y_1+\nu(u)(u-y_1)}+L_x\widehat{\phi}_n^{(c)}L_x\widehat{\phi}_{i-1}^{(c^{\prime})}L_{y_1}-\nu(u)L_x\widehat{\phi}_n^{(c)}L_x\widehat{\phi}_{i-1}^{(c^{\prime})}L_{y_1} \\
\displaystyle &{}& \quad -(-1)^{\nu(u)}\nu(u)L_x\widehat{\phi}_n^{(c)}(L_x\widehat{\phi}_{i-1}^{(c^{\prime})}L_{y_1+\nu(u)(u-y_1)}+\nu(u)L_u\widehat{\phi}_{i-1}^{(c^{\prime})}L_{u-y_1}) \\
\displaystyle &=& L_x\bigl([\widehat{\phi}_n^{(c)},L_{y_1}]\phi_{i-1}^{(c^{\prime})}-[\phi_{i-1}^{(c^{\prime})},L_{y_1}]\widehat{\phi}_n^{(c)}\bigl)L_{y_1+\nu(u)(u-y_1)} \\
\displaystyle &{}& \quad +L_x\bigl([\widehat{\phi}_n^{(c)},L_x]\phi_{i-1}^{(c^{\prime})}-[\phi_{i-1}^{(c^{\prime})},L_x]\widehat{\phi}_n^{(c)}\bigl)L_{y_1} \\
\displaystyle &{}& \quad -\nu(u)L_x\bigl([\widehat{\phi}_n^{(c)},L_x]\phi_{i-1}^{(c^{\prime})}-[\phi_{i-1}^{(c^{\prime})},L_x]\widehat{\phi}_n^{(c)}\bigl)L_{y_1} \\
\displaystyle &{}& \quad +\nu(u)L_x\bigl([\widehat{\phi}_n^{(c)},L_x]\phi_{i-1}^{(c^{\prime})}-[\phi_{i-1}^{(c^{\prime})},L_x]\widehat{\phi}_n^{(c)}\bigl)L_{y_1+\nu(u)(u-y_1)} \\
\displaystyle &{}& \quad +\nu(u)L_x\bigl([\widehat{\phi}_n^{(c)},L_u]\phi_{i-1}^{(c^{\prime})}-[\phi_{i-1}^{(c^{\prime})},L_u]\widehat{\phi}_n^{(c)}\bigl)L_{u-y_1} \\
\displaystyle &=& \nu(u)L_x\bigl([\widehat{\phi}_n^{(c)},L_u]\phi_{i-1}^{(c^{\prime})}-[\phi_{i-1}^{(c^{\prime})},L_u]\widehat{\phi}_n^{(c)}\bigl)L_{u-y_1}.
\end{eqnarray*}
On the other hand, 
\begin{eqnarray*}
\displaystyle Q &=& \widehat{\psi}_{n+1}^{(c)}(u)\widehat{\phi}_{i-1}^{(c^{\prime})}L_{u-y_1}+L_u\widehat{\phi}_{i-1}^{(c^{\prime})}\widehat{\psi}_{n+1}^{(c)}(u-y_1) \\
\displaystyle &{}& \quad -\widehat{\psi}_i^{(c^{\prime})}(u)\widehat{\phi}_n^{(c)}L_{u-y_1}-L_u\widehat{\phi}_n^{(c)}\widehat{\psi}_i^{(c^{\prime})}(u-y_1) \\
\displaystyle &=& (-1)^{\nu(u)}\bigl(L_x\widehat{\phi}_n^{(c)}L_{y_1+\nu(u)(u-y_1)}+\nu(u)L_u\widehat{\phi}_n^{(c)}L_{u-y_1}\bigl)\widehat{\phi}_{i-1}^{(c^{\prime})}L_{u-y_1} \\
\displaystyle &{}& \quad +(-1)^{\nu(u)}L_u\widehat{\phi}_{i-1}^{(c^{\prime})}\bigl(L_x\widehat{\phi}_n^{(c)}L_{y_1+\nu(u)(u-y_1)}+\nu(u)L_u\widehat{\phi}_n^{(c)}L_{u-y_1}\bigl) \\
\displaystyle &{}& \quad +L_u\widehat{\phi}_{i-1}^{(c^{\prime})}L_x\widehat{\phi}_n^{(c)}L_{y_1} \\
\displaystyle &{}& \quad -(-1)^{\nu(u)}(L_x\widehat{\phi}_{i-1}^{(c^{\prime})}L_{y_1+\nu(u)(u-y_1)}+\nu(u)L_u\widehat{\phi}_{i-1}^{(c^{\prime})}L_{u-y_1})\widehat{\phi}_n^{(c)}L_{u-y_1} \\
\displaystyle &{}& \quad -(-1)^{\nu(u)}L_u\widehat{\phi}_n^{(c)}(L_x\widehat{\phi}_{i-1}^{(c^{\prime})}L_{y_1+\nu(u)(u-y_1)}+\nu(u)L_u\widehat{\phi}_{i-1}^{(c^{\prime})}L_{u-y_1}) \\
\displaystyle &{}& \quad -L_u\widehat{\phi}_n^{(c)}L_x\widehat{\phi}_{i-1}^{(c^{\prime})}L_{y_1} \\
\displaystyle &=& (-1)^{\nu(u)}L_x\bigl([\widehat{\phi}_n^{(c)},L_{y_1+\nu(u)(u-y_1)}]\widehat{\phi}_{i-1}^{(c^{\prime})}-[\widehat{\phi}_{i-1}^{(c^{\prime})},L_{y_1+\nu(u)(u-y_1)}]\widehat{\phi}_n^{(c)}\bigl)L_{u-y_1} \\
\displaystyle &{}& \quad -\nu(u)L_u\bigl([\widehat{\phi}_n^{(c)},L_{u-y_1}]\widehat{\phi}_{i-1}^{(c^{\prime})}-[\widehat{\phi}_{i-1}^{(c^{\prime})},L_{u-y_1}]\widehat{\phi}_n^{(c)}\bigl)L_{u-y_1} \\
\displaystyle &{}& \quad +(-1)^{\nu(u)}L_u\bigl([\widehat{\phi}_n^{(c)},L_x]\widehat{\phi}_{i-1}^{(c^{\prime})}-[\widehat{\phi}_{i-1}^{(c^{\prime})},L_x]\widehat{\phi}_n^{(c)}\bigl)L_{y_1+\nu(u)(u-y_1)} \\
\displaystyle &{}& \quad -\nu(u)L_u\bigl([\widehat{\phi}_n^{(c)},L_u]\widehat{\phi}_{i-1}^{(c^{\prime})}-[\widehat{\phi}_{i-1}^{(c^{\prime})},L_u]\widehat{\phi}_n^{(c)}\bigl)L_{u-y_1} \\
\displaystyle &{}& \quad +L_u\bigl([\widehat{\phi}_n^{(c)},L_x]\widehat{\phi}_{i-1}^{(c^{\prime})}-[\widehat{\phi}_{i-1}^{(c^{\prime})},L_x]\widehat{\phi}_n^{(c)}\bigl)L_{y_1}.
\end{eqnarray*}
Thus, we obtain $\nu(u)Q=-\nu(u)L_x\bigl([\widehat{\phi}_n^{(c)},L_u]\phi_{i-1}^{(c^{\prime})}-[\phi_{i-1}^{(c^{\prime})},L_u]\widehat{\phi}_n^{(c)}\bigl)L_{u-y_1}=-P$, and (\ref{eq15}) is proven. 
\end{proof}

As an application, we obtain the following property.
\begin{prop}\label{prop4-6}
{\it We have
$$ \widehat{\partial}_n^{(c)}L_{x+\delta(s)y_s}(\mathcal{A}_{\mu_r}^1)\subset L_{x+\delta(s)y_s}(\mathcal{A}_{\mu_r}^1). $$
for any integer $n\ge 1$, any $c \in \mathbb{Q}$, and any $s\in\mu_r$. }
\end{prop}
\begin{proof}
The proof is given by induction on $n$. Based on the definition of $\widehat{\partial}_1^{(c)}=\partial_1$, the proposition holds for $n=1$. Suppose the assertion is true until $n-1$. We prove the case of $n$ by induction on the degrees of words. 

Based on $(\mathrm{A}_n)$ in the proof of Proposition \ref{prop4-5}, for $u\in\mathbb{Q}\cdot x+\sum_{s\in\mu_r}\mathbb{Q}\cdot y_s$, we have
\begin{eqnarray*}
\displaystyle \widehat{\partial}_n^{(c)}(u)&=&(-1)^{\nu(u)}(L_x\widehat{\phi}_{n-1}^{(c)}(y_1+\nu(u)(u-y_1))+\nu(u)L_u\widehat{\phi}_{n-1}^{(c)}(u-y_1)) \\
&=&
\left\{
\begin{array}{ll}
\displaystyle L_x\widehat{\phi}_{n-1}^{(c)}(y_1) & u=x, \\
\displaystyle -L_x\widehat{\phi}_{n-1}^{(c)}(y_1) & u=y_1, \\
\displaystyle -L_x\widehat{\phi}_{n-1}^{(c)}(y_s)-L_{y_s}\widehat{\phi}_{n-1}^{(c)}(y_s-y_1) & u=y_s,s\neq 1.
\end{array}
\right.
\end{eqnarray*}
By $(\beta_n)$ in the proof of Proposition \ref{prop4-5} and the induction hypothesis, we have $\widehat{\partial}_n^{(c)}(u)\in L_{x+\delta(s)y_s}(\mathcal{A}_{\mu_r}^1)$, and the assertion is proven for degree-$1$ words. 

Suppose this is proven until degree $<d$ and let $\deg(uw)=d$ with $u\in\mathbb{Q}\cdot x+\sum_{s\in\mu_r}\mathbb{Q}\cdot y_s$. Based on $(\mathrm{A}_n)$, we have
\begin{eqnarray*}
\displaystyle \widehat{\partial}_n^{(c)}(uw)&=&L_u\widehat{\partial}_n^{(c)}(w)+(-1)^{\nu(u)}(L_x\widehat{\phi}_{n-1}^{(c)}L_{y_1+\nu(u)(u-y_1)}+\nu(u)L_u\widehat{\phi}_{n-1}^{(c)}L_{u-y_1})(w) \\
&=&
\left\{
\begin{array}{ll}
\displaystyle L_x\widehat{\partial}_n^{(c)}(w)+L_x\widehat{\phi}_{n-1}^{(c)}L_{y_1}(w) & u=x, \\
\displaystyle L_{y_1}\widehat{\partial}_n^{(c)}(w)-L_x\widehat{\phi}_{n-1}^{(c)}L_{y_1}(w) & u=y_1, \\
\displaystyle L_{y_s}\widehat{\partial}_n^{(c)}(w)-L_x\widehat{\phi}_{n-1}^{(c)}L_{y_s}(w)-L_{y_s}\widehat{\phi}_{n-1}^{(c)}L_{y_s-y_1}(w) & u=y_s,s\neq 1.
\end{array}
\right.
\end{eqnarray*}
Accordingly, we see
$$ \displaystyle \widehat{\partial}_n^{(c)}(xw)\in x\mathcal{A}_{\mu_r},\widehat{\partial}_n^{(c)}(y_sw)\in x\mathcal{A}_{\mu_r}+ y_s\mathcal{A}_{\mu_r}~(s\neq 1). $$
Moreover, if $uw=w^{\prime}y_s$ $(s\in\mu_r)$, we obtain
$$ \displaystyle \widehat{\partial}_n^{(c)}(w^{\prime}y_s)\in\sum_{s\in\mu_r}\mathcal{A}_{\mu_r}y_s $$
based on $(\beta_n)$ and the induction hypothesis. Combining these statements, we obtain the assertion for degree-$d$ words. 
\end{proof}

\subsection{Extended Derivation Relation for MLV's}\label{subsec4-2}

\noindent In this subsection, we show the extended derivation relation for MLV's by reducing the relation to Corollary \ref{cor3-13}. 

Denote by $\mathcal{A}_{\mu_r,n}^1$ the weight-$n$ homogenous part of $\mathcal{A}_{\mu_r}^1$. Recall $z_{k,s}=x^{k-1}y_s$ for $k\ge 1$ and $s\in\mu_r$ as defined in $\S \ref{subsec3-3}$. Let $\mathfrak{W}$ be the $\mathbb{Q}$-vector space generated by $\{\mathcal{H}_w|w\in\mathcal{A}_{\{1\}}^1\}$, where $\mathcal{H}_w$ denotes the $\mathbb{Q}$-linear operator given by $\mathcal{H}_w(w^{\prime})=w\ast w^{\prime}$, and let $\mathfrak{W}_n$ be the vector subspace of $\mathfrak{W}$ generated by $\{\mathcal{H}_w|w\in\mathcal{A}_{\{1\},n}^1\}$. Let $\mathfrak{W}^{\prime}$ be the $\mathbb{Q}$-vector space generated by $\{L_{z_{k,1}}\mathcal{H}_w|1\le k,~w\in\mathcal{A}_{\{1\}}^1\}$, and let $\mathfrak{W}^{\prime}_n$ be the vector subspace of $\mathfrak{W}^{\prime}$ generated by $\{L_{z_{k,1}}\mathcal{H}_w|1\le k\le n,~w\in\mathcal{A}_{\{1\},n-k}^1\}$. We define a $\mathbb{Q}$-linear map $ \lambda : \mathfrak{W}^{\prime}\to\mathfrak{W} $ by $\lambda(L_{z_{k,1}}\mathcal{H}_w)=\mathcal{H}_{z_{k,1}w}$. 
\begin{rmk}\label{rmk4-7}
The map $\lambda$ is well-defined for the reasons described in the case of MZV's (see \cite{T}). 
\end{rmk}

\begin{lem}\label{lem4-8}
{\it For any $X\in\mathfrak{W}^{\prime}$, any $k\ge 1$, and any $s\in\mu_r$, we have 
$$ [\lambda(X),L_{z_{k,s}}]=XL_{z_{k,s}}+M_sL_{x^k}X. $$}
\end{lem}
\begin{proof}
It is sufficient to show the case of $X=L_{z_{k,1}}\mathcal{H}_w$, but it follows straightforward from the harmonic product rule:
\begin{equation}
[\mathcal{H}_{z_{k,1}w},L_{z_{l,s}}]=L_{z_{k,1}}\mathcal{H}_wL_{z_{l,s}}+M_sL_{x^l}L_{z_{k,1}}\mathcal{H}_w~~(s\in\mu_r). \label{eq18}
\end{equation}

\end{proof}
\begin{lem}\label{lem4-9}
{\it For any $k,l\ge 1$ and any $s\in\mu_r$, we have
$$ (\lambda -1)(\mathfrak{W}^{\prime}_k)L_{z_{l,s}}\subset M_s\cdot \mathfrak{W}^{\prime}_{k+l}. $$}
\end{lem}
\begin{proof}
This lemma is given by interpreting (\ref{eq18}) as
$$ \mathcal{H}_{z_{k,1}w}L_{z_{l,s}}-L_{z_{k,1}}\mathcal{H}_wL_{z_{l,s}}=M_s(L_{z_{l,1}}\mathcal{H}_{z_{k,1}w}+L_{x^l}L_{z_{k,1}}\mathcal{H}_w). $$
\end{proof}
\begin{lem}\label{lem4-10}
{\it For any $k,l\ge 1$ and any $s\in\mu_r$, we have
$$ (\lambda -1)(\mathfrak{W}^{\prime}_k)\cdot(\lambda -1)(\mathfrak{W}^{\prime}_l)\subset(\lambda -1)(\mathfrak{W}^{\prime}_{k+l}). $$}
\end{lem}
\begin{proof}
Let $d$ and $d^{\prime}$ be the weights of words $w$ and $w^{\prime}$, respectively. We need only show that $ (\lambda -1)(L_{z_k}\mathcal{H}_w)\cdot(\lambda -1)(L_{z_l}\mathcal{H}_{w^{\prime}})\in (\lambda -1)(\mathfrak{W}^{\prime}_{k+l+d+d^{\prime}}). $
\[
\begin{array}{ccl}
\text{LHS} &=& (\mathcal{H}_{z_{k,1}w}-L_{z_{k,1}}\mathcal{H}_w)(\mathcal{H}_{z_{l,1}w^{\prime}}-L_{z_{l,1}}\mathcal{H}_{w^{\prime}}) \\
&=& \mathcal{H}_{z_{k,1}w\ast z_{l,1}w^{\prime}}-\mathcal{H}_{z_{k,1}w}L_{z_{l,1}}\mathcal{H}_{w^{\prime}}-L_{z_{k,1}}\mathcal{H}_{w\ast z_{l,1}w^{\prime}}+L_{z_{k,1}}\mathcal{H}_wL_{z_{l,1}}\mathcal{H}_{w^{\prime}} \\
&=& \mathcal{H}_{z_{k,1}(w\ast z_{l,1}w^{\prime})+z_{l,1}(z_{k,1}w\ast w^{\prime})+z_{k+l,1}(w\ast w^{\prime})}-(L_{z_{k,1}}\mathcal{H}_wL_{z_{l,1}} \\
&{}& \quad +L_{z_{l,1}}\mathcal{H}_{z_{k,1}w}+L_{z_{k+l,1}}\mathcal{H}_w)\mathcal{H}_{w^{\prime}}-L_{z_{k,1}}\mathcal{H}_{w\ast z_{l,1}w^{\prime}}+L_{z_{k,1}}\mathcal{H}_wL_{z_{l,1}}\mathcal{H}_{w^{\prime}} \\
&=& \mathcal{H}_{z_{k,1}(w\ast z_{l,1}w^{\prime})}-L_{z_{k,1}}\mathcal{H}_{w\ast z_{l,1}w^{\prime}}+\mathcal{H}_{z_{l,1}(z_{k,1}w\ast w^{\prime})}-L_{z_{l,1}}\mathcal{H}_{z_{k,1}w\ast w^{\prime}} \\
&{}& \quad +\mathcal{H}_{z_{k+l,1}(w\ast w^{\prime})}-L_{z_{k+l,1}}\mathcal{H}_{w\ast w^{\prime}} \\
&=& (\lambda -1)(L_{z_{k,1}}\mathcal{H}_{w\ast z_{l,1}w^{\prime}}+L_{z_{l,1}}\mathcal{H}_{z_{k,1}w\ast w^{\prime}}+L_{z_{k+l,1}}\mathcal{H}_{w\ast w^{\prime}}). \\
&\in& \mathrm{RHS}.
\end{array}
\]
Hence, the lemma is proven. 
\end{proof}
\begin{lem}\label{lem4-11}
{\it For any $X\in\mathfrak{W}^{\prime}$, we have $\lambda(X)(1)=X(1)$.}
\end{lem}
\begin{proof}
$$(\lambda-1)(L_{z_{k,1}}\mathcal{H}_w)(1)=\mathcal{H}_{z_{k,1}w}(1)-L_{z_{k,1}}\mathcal{H}_w(1)=z_{k,1}w-z_{k,1}w=0.$$
\end{proof}
\begin{lem}\label{lem4-12}
{\it Let $X\in\mathfrak{W}$. If $X(1)=0$ and $[X,L_{z_{k,s}}]=0$ for any $k\ge 1$ and any $s\in\mu_r$, we have $X=0$ on $\mathcal{A}_{\mu_r}^1$.}
\end{lem}
\begin{proof}
If $[X,L_{z_{k,s}}]=0$ for any $k\ge 1$ and any $s\in\mu_r$,
$$ X(z_{k_1,s_1}\cdots z_{k_n,s_n})=z_{k_1,s_1}X(z_{k_2,s_2}\cdots z_{k_n,s_n}) = \cdots = z_{k_1,s_1}\cdots z_{k_n,s_n}X(1)=0. $$
\end{proof}

Now, let $\sigma_s$ ($s\in\mu_r$) denote $\varphi\mathcal{I}M_s$. The following proposition is the key to connect the extended derivation operator $\widehat{\partial}_n^{(c)}$ with the harmonic product operator $\mathcal{H}_w$. 
\begin{prop}\label{prop4-13}
{\it Let $n\ge 1$, $c\in\mathbb{Q}$ and $s\in\mu_r$. Then, the following two statements, $(\mathrm{C}_n)$ and $(\mathrm{D}_n)$ hold. }
\begin{itemize}
\item[$(\mathrm{C}_n)$]$\sigma_s^{-1}\widehat{\phi}_{n-1}^{(c)}L_{y_1-\delta(s)y_s}\sigma_s\in\mathfrak{W}_n^{\prime}, $
\item[$(\mathrm{D}_n)$]$\sigma_s^{-1}L_{x+\delta(s)y_s}^{-1}\widehat{\partial}_n^{(c)}L_{x+\delta(s)y_s}\sigma_s =\lambda(\sigma_s^{-1}\widehat{\phi}_{n-1}^{(c)}L_{y_1-\delta(s)y_s}\sigma_s)~(\in\mathfrak{W}_n)~~~\mathit{on}~\mathcal{A}_{\mu_r}^1. $
\end{itemize}
\end{prop}
\noindent Note that, by Proposition \ref{prop4-6}, the expression $L_{x+\delta(s)y_s}^{-1}$ in $(\mathrm{D}_n)$ has the well-defined meaning. 

\begin{proof}
We prove inductively that $(\mathrm{C}_1)\Rightarrow(\mathrm{D}_1)\Rightarrow(\mathrm{C}_2)\Rightarrow(\mathrm{D}_2)\Rightarrow(\mathrm{C}_3)\Rightarrow\cdots$. Since
$$ \sigma_s^{-1}\widehat{\phi}_0^{(c)}L_{y_1-\delta(s)y_s}\sigma_s=-\sigma_s^{-1}\varphi L_{y_s}\mathcal{I}M_s=-\sigma_s^{-1}\varphi\mathcal{I}M_sL_{y_1}=-L_{y_1}, $$
the assertion $(\mathrm{C}_1)$ holds. Suppose that the above statement is true until $(\mathrm{C}_n)$. Note that
\begin{equation}
L_{x+\delta(s)y_s}^{-1}\widehat{\partial}_n^{(c)}L_{x+\delta(s)y_s}=\widehat{\partial}_n^{(c)}+\widehat{\phi}_{n-1}^{(c)}L_{y_1-\delta(s)y_s} \label{eq16}
\end{equation}
based on $(\mathrm{A}_n)$ in the proof of Proposition \ref{prop4-5} and Corollary \ref{cor4-4} (ii). In addition, we note that $(\mathrm{C}_n)$ implies the operator $\sigma_s^{-1}\widehat{\phi}_{n-1}^{(c)}L_{y_1-\delta(s)y_s}\sigma_s$ is independent of $s\in\mu_r$. 

First, we prove that the operator $\sigma_s^{-1}\widehat{\partial}_n^{(c)}\sigma_s$ is also independent of $s\in\mu_r$ by induction on the $y_s$-degree of a word. If $w=x^l$ ($l\ge 0$), we have
$$ \sigma_s^{-1}\widehat{\partial}_n^{(c)}\sigma_s(w)=\sigma_s^{-1}\widehat{\partial}_n^{(c)}(z^l)=0 $$
and hence the claim holds. Suppose that $\sigma_s^{-1}\widehat{\partial}_n^{(c)}\sigma_s(w^{\prime})$ ($w^{\prime}$ is a word of $\mathcal{A}_{\mu_r}$) that is independent of $s\in\mu_r$, and let $w=z_{l,t}w^{\prime}$ ($l\ge 1,t\in\mu_r$). We find that
\begin{eqnarray*}
\sigma(w)&=&\varphi\mathcal{I}M_sL_{z_{l,t}}(w^{\prime}) \\
&=&\varphi L_{z_{l,st}}\mathcal{I}M_{st}(w^{\prime}) \\
&=&L_{z^{l+1}}L_{\delta(st)y_{st}-y_1}\sigma_{st}(w^{\prime})
\end{eqnarray*}
Then, assuming $(\mathrm{A}_n)$ and Corollary \ref{cor4-4}, we have
\begin{eqnarray*}
&{}&\sigma_s^{-1}\widehat{\partial}_n^{(c)}\sigma_s(w) \\
&=&\sigma_s^{-1}\widehat{\partial}_n^{(c)}L_{z^{l-1}}L_{\delta(st)y_{st}-y_1}\sigma_{st}(w^{\prime}) \\
&=&\sigma_s^{-1}L_{z^{l-1}}(L_{\delta(st)y_{st}-y_1}\widehat{\partial}_n^{(c)}+\widehat{\psi}_n^{(c)}(\delta(st)y_{st}-y_1))\sigma_{st}(w^{\prime}) \\
&=&L_{z_{l,t}}\sigma_{st}^{-1}\widehat{\partial}_n^{(c)}\sigma_{st}(w^{\prime})+(M_tL_{x^l}+L_{z_{l,t}})\sigma_{st}^{-1}\widehat{\phi}_{n-1}^{(c)}L_{y_1-\delta(st)y_{st}}\sigma_{st}(w^{\prime}).
\end{eqnarray*}
Hence, $\sigma_s^{-1}\widehat{\partial}_n^{(c)}\sigma_s(w)$ is independent of $s\in\mu_r$. 

Now, we see that
\begin{eqnarray*}
&{}&[\sigma_s^{-1}L_{x+\delta(s)y_s}^{-1}\widehat{\partial}_n^{(c)}L_{x+\delta(s)y_s}\sigma_s,L_{z_{l,t}}] \\
&=&\sigma_s^{-1}(\widehat{\partial}_n^{(c)}+\widehat{\phi}_{n-1}^{(c)}L_{y_1-\delta(s)y_s})\sigma_sL_{z_{l,t}}-L_{z_{l,t}}\sigma_s^{-1}(\widehat{\partial}_n^{(c)}+\widehat{\phi}_{n-1}^{(c)}L_{y_1-\delta(s)y_s})\sigma_s.
\end{eqnarray*}
Since
\begin{eqnarray*}
&{}&\sigma_s^{-1}\widehat{\partial}_n^{(c)}\sigma_sL_{z_{l,t}} \\
&=&\sigma_s^{-1}\widehat{\partial}_n^{(c)}L_{z^{l-1}}L_{\delta(st)y_{st}-y_1}\sigma_s \\
&=&\sigma_s^{-1}L_{z^{l-1}}(L_{\delta(st)y_{st}-y_1}\widehat{\partial}_n^{(c)}+\widehat{\psi}_n^{(c)}(\delta(st)y_{st}-y_1))\sigma_{st} \\
&=&L_{z_{l,t}}\sigma_s^{-1}\widehat{\partial}_n^{(c)}\sigma_{st}+(M_tL_{x^l}+L_{z_{l,t}})\sigma_{st}^{-1}\widehat{\phi}_{n-1}^{(c)}L_{y_1-\delta(st)y_{st}}\sigma_{st} \\
&=&L_{z_{l,t}}\sigma_s^{-1}\widehat{\partial}_n^{(c)}\sigma_s+(M_tL_{x^l}+L_{z_{l,t}})\sigma_s^{-1}\widehat{\phi}_{n-1}^{(c)}L_{y_1-\delta(s)y_s}\sigma_s,
\end{eqnarray*}
we have
\begin{eqnarray*}
&{}&[\sigma_s^{-1}L_{x+\delta(s)y_s}^{-1}\widehat{\partial}_n^{(c)}L_{x+\delta(s)y_s}\sigma_s,L_{z_{l,t}}] \\
&=&\sigma_s^{-1}\widehat{\phi}_{n-1}^{(c)}L_{y_1-\delta(s)y_s}\sigma_sL_{z_{l,t}}+M_tL_{x^l}\sigma_s^{-1}\widehat{\phi}_{n-1}^{(c)}L_{y_1-\delta(s)y_s}\sigma_s.
\end{eqnarray*}
According to Lemma \ref{lem4-8}, this is equivalent to $[\lambda(\sigma_s^{-1}\widehat{\phi}_{n-1}^{(c)}L_{y_1-\delta(s)y_s}\sigma_s),L_{z_{l,t}}]$. In addition, we obtain
\begin{eqnarray*}
&{}&\sigma_s^{-1}L_{x+\delta(s)y_s}^{-1}\widehat{\partial}_n^{(c)}L_{x+\delta(s)y_s}\sigma_s(1) \\
&=&\sigma_s^{-1}(\widehat{\partial}_n^{(c)}+\widehat{\phi}_{n-1}^{(c)}L_{y_1-\delta(s)y_s})(1) \\
&=& \sigma_s^{-1}\widehat{\phi}_{n-1}^{(c)}L_{y_1-\delta(s)y_s}\sigma_s(1).
\end{eqnarray*}
Hence, by Lemma \ref{lem4-12}, we have $(\mathrm{D}_n)$. 

On the other hand, suppose that $(\mathrm{D}_n)$ is proven. Using (\ref{eq16}) and $(\mathrm{D}_n)$, we find
\begin{equation}
\sigma_s^{-1}\widehat{\partial}_n^{(c)}\sigma_s=(\lambda-1)(\sigma_s^{-1}\widehat{\phi}_{n-1}^{(c)}L_{y_1-\delta(s)y_s}\sigma_s). \label{eq17}
\end{equation}
In addition, based on $(\beta_n)$, we have the expression
$$ \displaystyle \widehat{\phi}_n^{(c)}=\sum_{i=0}^n\widehat{f}_i^{(c)}L_{z^{n-i}}~~(\widehat{f}_i^{(c)}\in\mathbb{Q}[\widehat{\partial}_1^{(c)},\ldots,\widehat{\partial}_n^{(c)}]_{(i)}). $$
Hence,
\begin{eqnarray*}
\displaystyle \sigma_s^{-1}\widehat{\phi}_{n-1}^{(c)}L_{y_1-\delta(s)y_s}\sigma_s&=&\sum_{i=0}^n\sigma_s^{-1}\widehat{f}_i^{(c)}L_{z^{n-i}}L_{y_1-\delta(s)y_s}\sigma_s \\
\displaystyle &=&-\sum_{i=0}^n\sigma_s^{-1}\widehat{f}_i^{(c)}\sigma_sL_{z_{n-i+1,1}}.
\end{eqnarray*}
By Lemma \ref{lem4-10} and (\ref{eq17}), this is an element of $\sum_{i=1}^n(\lambda-1)(\mathfrak{W}_i^{\prime})L_{z_{n-i+1,1}}+\mathbb{Q}\cdot L_{z_{n,1}}$. Then, by Lemma \ref{lem4-9}, this is a subset of $\mathfrak{W}_{n+1}^{\prime}$. Thus, $(\mathrm{C}_{n+1})$ is proven. 
\end{proof}

\begin{thm}\label{thm4-14}
{\it For any $n\ge 1$ and any $c\in\mathbb{Q}$, we have $\widehat{\partial}_n^{(c)}(\mathcal{A}_{\mu_r}^0)\subset \mathrm{Ker}(\mathcal{L}^{\mathcyb{sh}}).$ }
\end{thm}

\begin{proof}
By Proposition \ref{prop4-13} $(\mathrm{D}_n)$, there exists $w\in\mathcal{A}_{\{1\},n}^1$ such that
$$ \widehat{\partial}_n^{(c)}L_{x+\delta(s)y_s}\sigma_s=L_{x+\delta(s)y_s}\sigma_s\mathcal{H}_w $$
on $\mathcal{A}_{\mu_r}^1$. ($w=L_{x+\delta(s)y_s}^{-1}\sigma_s^{-1}\widehat{\partial}_n^{(c)}(x+\delta(s)y_s)=L_x^{-1}\varphi\widehat{\partial}_n^{(c)}(x)$.) Since
$$ L_{x+\delta(s)y_s}\sigma_s(\mathcal{A}_{\mu_r,>0}^1)=\mathcal{A}_{\mu_r,>0}^0, $$
we find that
$$ \widehat{\partial}_n^{(c)}(\mathcal{A}_{\mu_r,>0}^0)=\widehat{\partial}_n^{(c)}L_{x+\delta(s)y_s}\sigma_s(\mathcal{A}_{\mu_r,>0}^1)=L_{x+\delta(s)y_s}\sigma_s\mathcal{H}_w(\mathcal{A}_{\mu_r,>0}^1). $$
This is a subset of $L_{x+\delta(s)y_s}\varphi\mathcal{I}M_s(\mathcal{A}_{\mu_r,>0}^1\ast\mathcal{A}_{\{1\},>0}^1)$. By Corollary \ref{cor3-13}, we obtain 
$$ \widehat{\partial}_n^{(c)}(\mathcal{A}_{\mu_r,>0}^0)\subset\mathrm{Ker}\mathcal{L}^{\mathcyb{sh}} $$
and hence the theorem because $\widehat{\partial}_n^{(c)}(\mathbb{Q})=\{0\}\subset\mathrm{Ker}\mathcal{L}^{\mathcyb{sh}}$.
\end{proof}

Next, we study the properties of the alternative operator $\partial_n^{(c)}$, which Kaneko devised by modeling a Hopf algebra developed by A. Connes and H. Moscovici (see \cite{CM} for details of the structure), and thereby find that $\partial_n^{(c)}$ also induces relations of MLV's. The commutativity property of $\partial_n^{(c)}$ was investigated by Wakiyama \cite{W} in advance. 
\begin{defn}\label{defn4-15}
{\it Let $c$ be a rational number, and let $H$ be the derivation on $\mathcal{A}_{\mu_r}$ defined in Definition \ref{defn4-1}. For each integer $n \ge 1$, we define a $\mathbb{Q}$-linear map $\partial_n^{(c)}$ from $\mathcal{A}_{\mu_r}$ to $\mathcal{A}_{\mu_r}$ by}
$$ \displaystyle \partial_n^{(c)}=\frac{1}{(n-1)!}\mathrm{ad}(\theta^{(c)})^{n-1}(\partial_1), $$
{\it where $\theta^{(c)}$ is the $\mathbb{Q}$-linear map defined by $\theta^{(c)}(u)=\theta (u)$ ($u=x$ or $y_s$) and the rule}
$$ \theta^{(c)}(ww^{\prime})=\theta^{(c)}(w)w^{\prime}+w\theta^{(c)}(w^{\prime})+c\partial_1(w)H(w^{\prime}) $$
{\it for any $w,w^{\prime} \in \mathcal{A}_{\mu_r}$.}
\end{defn}
\noindent The only difference between $\widehat{\theta}^{(c)}$ and $\theta^{(c)}$ is the order of $H$ and $\partial_1$ appearing in the right-hand side of their recursive rules. 
\begin{lem}\label{lem4-16}
{\it For any rational number $c$, we have
$$ \theta^{(c)}=\widehat{\theta}^{(-c)}+c\partial_1(H-1). $$}
\end{lem}

\begin{proof}
We immediately see that the images of generators $x$ and $y_s$ $(s\in\mu_r)$ of both sides coincide. For $w,w^{\prime}\in\mathcal{A}_{\mu_r}$, 
\begin{eqnarray*}
&{}&(\widehat{\theta}^{(-c)}+c\partial_1(H-1))(ww^{\prime}) \\
&=&\widehat{\theta}^{(-c)}(w)w^{\prime}+w\widehat{\theta}^{(-c)}(w^{\prime})-cH(w)\partial_1(w^{\prime})+c(\partial_1H(w)w^{\prime} \\
&{}&\quad +H(w)\partial_1(w^{\prime})+\partial_1(w)H(w^{\prime})+w\partial_1H(w^{\prime}))-c(\partial_1(w)w^{\prime}+w\partial_1(w^{\prime})) \\
&=&(\widehat{\theta}^{(-c)}+c\partial_1(H-1))(w)w^{\prime}+w(\widehat{\theta}^{(-c)}+c\partial_1(H-1))(w^{\prime})+c\partial_1(w)H(w^{\prime}).
\end{eqnarray*}
Hence, the recursive rules also coincide, and the Lemma is proven. 
\end{proof}

\begin{prop}\label{prop4-17}
{\it For any positive integer $n$ and any rational number $c$, we have
$$ \partial_n^{(c)}\in \mathbb{Q}[\widehat{\partial}_1^{(-c)},\ldots ,\widehat{\partial}_n^{(-c)}]_{(n)}. $$
}
\end{prop}

\begin{proof}
The proposition holds for $n=1$ because $\partial_1^{(c)}=\widehat{\partial}_1^{(-c)}=\partial_1$. Suppose that the proposition is proven for $n$. Using Lemma \ref{lem4-16} and Proposition \ref{prop4-5}, we obtain
\[
\begin{array}{ccl}
n\partial_{n+1}^{(c)} &=& [\theta^{(c)},\partial_n^{(c)}] \\
&=& [\widehat{\theta}^{(-c)}+c\partial_1(H-1),\partial_n^{(c)}] \\
&=& [\widehat{\theta}^{(-c)},\partial_n^{(c)}]+c(n-1)\partial_1\partial_n^{(c)}.\end{array}
\]
Hence, by induction the Proposition holds for $n+1$. 
\end{proof}

\begin{cor}\label{cor4-18}
{\it For any rational numbers $c, c^{\prime}$, and any positive integers $n, m$, we have
$$ [\partial_n^{(c)},\widehat{\partial}_m^{(c^{\prime})}]=0. $$ }
\end{cor}

\begin{proof}
Immediately from Proposition \ref{prop4-5} and \ref{prop4-17}. 
\end{proof}
\noindent By Proposition \ref{prop4-17} and Theorem \ref{thm4-14}, we have
\begin{equation}
\partial_n^{(c)}(\mathcal{A}_{\mu_r}^0)\subset\mathrm{Ker}\mathcal{L}^{\mathcyb{sh}}, \label{eq}
\end{equation}
which gives the same class as Theorem \ref{thm4-14}. 

\subsection{On the Derivation Relation}\label{subsec4-3}

\noindent Let $\widehat{\mathcal{A}}_{\mu_r}$ denote the completion of $\mathcal{A}_{\mu_r}$. We put $\widehat{\Delta}$ as
\begin{equation}
\displaystyle\widehat{\Delta}=\exp\biggl(\sum_{n=1}^{\infty}\frac{\partial_n}{n}\biggl) \label{eq13}
\end{equation}
which has been introduced in \cite{AK}. Then, $\widehat{\Delta}$ is an automorphism on $\widehat{\mathcal{A}}_{\mu_r}$ that satisfies
$$ \displaystyle \widehat{\Delta}(x)=x\frac{1}{1-y_1},~\widehat{\Delta}(y_s)=(x+y_s)\frac{1}{1-y_1+y_s} $$
for any $s\in\mu_r$. We find that
$$ \displaystyle \widehat{\Delta}(x+\delta(s)y_s)=(x+\delta(s)y_s)\frac{1}{1-y_1+\delta(s)y_s},~\widehat{\Delta}(z)=z, $$
or, in other words, 
$$ \displaystyle \widehat{\Delta}L_{x+\delta(s)y_s}=L_{x+\delta(s)y_s}L_{\frac{1}{1-y_1+\delta(s)y_s}}\widehat{\Delta},~L_z\widehat{\Delta}=\widehat{\Delta}L_z. $$

Now, let $\Phi$ be
$$ \displaystyle \Phi(w)=(1+y_1)\biggl(\frac{1}{1+y_1}\ast w\biggl)\quad (w\in\widehat{\mathcal{A}}_{\mu_r}^1) $$
as in \cite{AK}, where $\widehat{\mathcal{A}}_{\mu_r}^1$ is the completion of $\mathcal{A}_{\mu_r}^1$. We see that $\Phi$ can be extended to an automorphism on $\widehat{\mathcal{A}}_{\mu_r}$ that satisfies
$$ \displaystyle \Phi(x)=x,~\Phi(x+y_s)=y_s\frac{1}{1+y_1}+x-xy_s\frac{1}{1+y_1} $$
and hence
$$ \Phi(y_s)=(1-x)y_s\frac{1}{1+y_1}. $$
\begin{lem}\label{lem4-19}
{\it Let $s\in\mu_r$. We have
$\displaystyle\varphi\mathcal{I}M_s\Phi=\widehat{\Delta}\varphi\mathcal{I}M_s$ on $\mathcal{A}_{\mu_r}$.}
\end{lem}

\begin{proof}
We need only show the equality of the images of general monomials in $\mathcal{A}_{\mu_r}$ on both sides. 
\begin{eqnarray*}
\displaystyle &{}& \widehat{\Delta}\varphi\mathcal{I}M_s(z_{k_1,s_1}\cdots z_{k_n,s_n}x^l) \\
\displaystyle &=& \widehat{\Delta}\varphi(z_{k_1,ss_1}\cdots z_{k_n,ss_1\cdots s_n}x^l) \\
\displaystyle &=& \widehat{\Delta}\bigl(z^{k_1-1}(\delta(ss_1)y_{ss_1}-y_1)\cdots z^{k_n-1}(\delta(ss_1\cdots s_n)y_{ss_1\cdots s_n}-y_1)z^l\bigl) \\
\displaystyle &=& z^{k_1-1}\biggl((x+\delta(ss_1)y_{ss_1})\frac{1}{1-y_1+\delta(ss_1)y_{ss_1}}-z\biggl) \\
\displaystyle &{}& \quad \cdots z^{k_n-1}\biggl((x+\delta(ss_1\cdots s_n)y_{ss_1\cdots s_n})\frac{1}{1-y_1+\delta(ss_1\cdots s_n)y_{ss_1\cdots s_n}}-z\biggl)z^l.
\end{eqnarray*}
On the other hand, 
\begin{eqnarray*}
\displaystyle &{}& \varphi\mathcal{I}M_s\Phi(z_{k_1,s_1}\cdots z_{k_n,s_n}x^l) \\
\displaystyle &=& \varphi\mathcal{I}M_s\biggl(x^{k_1-1}(1-x)y_{s_1}\frac{1}{1+y_1}\cdots x^{k_n-1}(1-x)y_{s_n}\frac{1}{1+y_1}x^l\biggl) \\
\displaystyle &=& \varphi\biggl(x^{k_1-1}(1-x)y_{ss_1}\frac{1}{1+y_{ss_1}}\cdots x^{k_n-1}(1-x)y_{ss_1\cdots s_n}\frac{1}{1+y_{ss_1\cdots s_n}}x^l\biggl) \\
\displaystyle &=& z^{k_1-1}(1-z)(\delta(ss_1)y_{ss_1}-y_1)\frac{1}{1-y_1+\delta(ss_1)y_{ss_1}} \\
\displaystyle &{}& \quad \cdots z^{k_n-1}(1-z)(\delta(ss_1\cdots s_n)y_{ss_1\cdots s_n}-y_1)\frac{1}{1-y_1+\delta(ss_1\cdots s_n)y_{ss_1\cdots s_n}}z^l.
\end{eqnarray*}
We can easily see that
$$ (x+\delta(s)y_s)\frac{1}{1-y_1+\delta(s)y_s}=z+(1-z)(\delta(s)y_s-y_1)\frac{1}{1-y_1+\delta(s)y_s}. $$
Hence, the lemma is proven.
\end{proof}

\begin{prop}\label{prop4-20}
{\it Let $s\in\mu_r$. Then, we have
$$ \displaystyle\mathcal{H}_{\frac{1}{1+y_1}}=M_{\frac{1}{s}}\mathcal{I}^{-1}\varphi L_{x+\delta(s)y_s}\widehat{\Delta}L_{x+\delta(s)y_s}\varphi\mathcal{I}M_s $$
on $\mathcal{A}_{\mu_r}^1.$}
\end{prop}

\begin{proof}
Using Lemma \ref{lem4-19}, we have
\begin{eqnarray*}
\displaystyle\mathrm{RHS} &=& M_{\frac{1}{s}}\mathcal{I}^{-1}\varphi L_{\frac{1}{1-y_1+\delta(f)y_f}}\widehat{\Delta}_{1-X}\varphi\mathcal{I}M_s \\
\displaystyle &=& M_{\frac{1}{s}}\mathcal{I}^{-1} L_{\frac{1}{1+y_s}}\varphi\widehat{\Delta}_{1-X}\varphi\mathcal{I}M_s \\
\displaystyle &=& L_{\frac{1}{1+y_1}}M_{\frac{1}{s}}\mathcal{I}^{-1} \varphi\widehat{\Delta}_{1-X}\varphi\mathcal{I}M_s \\
\displaystyle &=& L_{\frac{1}{1+y_1}}\Phi_{-y_1} = \mathrm{LHS}.
\end{eqnarray*}
\end{proof}
\noindent Based on equation (\ref{eq13}), Proposition \ref{prop4-20}, and Corollary \ref{cor3-13}, we again obtain the derivation relation
$$ \partial_n(\mathcal{A}_{\mu_r}^0)\subset\mathrm{Ker}\mathcal{L}^{\mathcyb{sh}}. $$


\section{Proofs of Lemmata}\label{sec5}

\begin{proof}[Proof of Lemma \ref{lem3-7}]
If (i) and (ii) are proven, then (i)${}^{\prime}$ and (ii)${}^{\prime}$ also hold by substituting $d_{\ast}^{-1}(w),d_{\ast}^{-1}(w^{\prime})$ in place of $w,w^{\prime}$ and by applying $d_{\ast}^{-1}$ to both sides. Therefore, we show (i) and (ii). 

(i)~By induction on the depth of a word. We denote the depth of a word $w$ by $\mathrm{dep}(w)$. When $\mathrm{dep}(w)=1$, set $w=z_{k,s},w^{\prime}=z_{l,t}v$. If $\mathrm{dep}(v)=0$, then
\begin{eqnarray*}
\mathrm{LHS} &=& \mathcal{I}^{-1}d_{\mathcyb{sh}}\mathcal{I}(z_{k,s}~\overline{\ast}~z_{l,t}) \\
&=& \mathcal{I}^{-1}d_{\mathcyb{sh}}\mathcal{I}(z_{k,s}z_{l,t}+z_{l,t}z_{k,s}-z_{k+l,st}) \\
&=& \mathcal{I}^{-1}d_{\mathcyb{sh}}(z_{k,s}z_{l,st}+z_{l,t}z_{k,st}-z_{k+l,st}) \\
&=& \mathcal{I}^{-1}((z_{k,s}+x^k)z_{l,st}+(z_{l,t}+x^l)z_{k,st}-z_{k+l,st}) \\
&=& \mathcal{I}^{-1}(z_{k,s}z_{l,st}+z_{l,t}z_{k,st}+z_{k+l,st})=\mathrm{RHS}.
\end{eqnarray*}
Hence, the assertion is proven. 

If $\mathrm{dep}(v)>0$, then
\begin{eqnarray*}
\mathrm{LHS} &=& \mathcal{I}^{-1}d_{\mathcyb{sh}}\mathcal{I}(z_{k,s}~\overline{\ast}~z_{l,t}v) \\
&=& \mathcal{I}^{-1}d_{\mathcyb{sh}}\mathcal{I}(z_{k,s}z_{l,t}v+z_{l,t}(z_{k,s}~\overline{\ast}~v)~-z_{k+l,st}v) \\
&=& \mathcal{I}^{-1}d_{\mathcyb{sh}}L_{z_{k,s}}N_{\frac{1}{s}}L_{z_{l,t}}N_{\frac{1}{t}}\mathcal{I}(v)+\mathcal{I}^{-1}d_{\mathcyb{sh}}L_{z_{l,t}}N_{\frac{1}{t}}\mathcal{I}(z_{k,s}~\overline{\ast}~v) \\
&{}& \quad -\mathcal{I}^{-1}d_{\mathcyb{sh}}L_{z_{k+l,st}}N_{\frac{1}{st}}\mathcal{I}(v) \\
&=& \mathcal{I}^{-1}L_{z_{k,s}+x^k}N_{\frac{1}{s}}L_{z_{l,t}+x^l}N_{\frac{1}{t}}d_{\mathcyb{sh}}\mathcal{I}(v)+\mathcal{I}^{-1}L_{z_{l,t}+x^l}N_{\frac{1}{t}}d_{\mathcyb{sh}}\mathcal{I}(z_{k,s}~\overline{\ast}~v) \\
&{}& \quad -\mathcal{I}^{-1}L_{z_{k+l,st}+x^{k+l}}N_{\frac{1}{st}}d_{\mathcyb{sh}}\mathcal{I}(v) \\
&=& L_{z_{k,s}}L_{z_{l,t}}\mathcal{I}^{-1}d_{\mathcyb{sh}}\mathcal{I}(v)+L_{z_{k,s}}L_{x^l}M_t\mathcal{I}^{-1}d_{\mathcyb{sh}}\mathcal{I}(v) \\
&{}& \quad +L_{x^k}L_{z_{l,st}}\mathcal{I}^{-1}d_{\mathcyb{sh}}\mathcal{I}(v)+L_{x^{k+l}}M_{st}\mathcal{I}^{-1}d_{\mathcyb{sh}}\mathcal{I}(v) \\
&{}& \quad +L_{z_{l,t}}(z_{k,s}\ast\mathcal{I}^{-1}d_{\mathcyb{sh}}\mathcal{I}(v))+L_{x^l}M_t(z_{k,s}~\overline{\ast}~\mathcal{I}^{-1}d_{\mathcyb{sh}}\mathcal{I}(v)) \\
&{}& \quad -L_{z_{k+l,st}}\mathcal{I}^{-1}d_{\mathcyb{sh}}\mathcal{I}(v)-L_{x^{k+l}}M_{st}\mathcal{I}^{-1}d_{\mathcyb{sh}}\mathcal{I}(v),
\end{eqnarray*}
\begin{eqnarray*}
\mathrm{RHS} &=& z_{k,s}\ast\mathcal{I}^{-1}d_{\mathcyb{sh}}\mathcal{I}(z_{l,t}v) \\
&=& z_{k,s}\ast(L_{z_{l,t}}\mathcal{I}^{-1}d_{\mathcyb{sh}}\mathcal{I}(v)+L_{x^l}M_t\mathcal{I}^{-1}d_{\mathcyb{sh}}\mathcal{I}(v)) \\
&=& L_{z_{k,s}}L_{z_{l,t}}\mathcal{I}^{-1}d_{\mathcyb{sh}}\mathcal{I}(v)+L_{z_{l,t}}(z_{k,s}\ast\mathcal{I}^{-1}d_{\mathcyb{sh}}\mathcal{I}(v)) \\
&{}& \quad +L_{z_{k+l,st}}\mathcal{I}^{-1}d_{\mathcyb{sh}}\mathcal{I}(v)+z_{k,s}\ast(L_{x^l}M_t\mathcal{I}^{-1}d_{\mathcyb{sh}}\mathcal{I}(v)).
\end{eqnarray*}
Hence, we must show that
$$ L_{z_{k,s}}L_{x^l}M_t(V)+L_{x^l}M_t(z_{k,s}\ast V)=L_{z_{k+l,st}}(V)+z_{k,s}\ast(L_{x^l}M_t(V)). $$
Setting $V=z_{k_1,s_1}\cdots z_{k_1,s_n}$, we have
\begin{eqnarray*}
\mathrm{LHS} &=& z_{k,s}z_{k_1+l,s_1t}z_{k_2,s_2}\cdots z_{k_1,s_n}+L_{x^l}M_t\bigl(z_{k,s}z_{k_1,s_1}\cdots z_{k_1,s_n} \\
&{}& \quad +z_{k_1,s_1}(z_{k,s}\ast z_{k_1,s_2}\cdots z_{k_1,s_n})+z_{k+k_1,ss_1}z_{k_1,s_2}\cdots z_{k_1,s_n}\bigl),
\end{eqnarray*}
\begin{eqnarray*}
\mathrm{RHS} &=& z_{k+l,st}z_{k_1,s_1}\cdots z_{k_1,s_n}+z_{k,s}\ast z_{k_1+l,s_1t}z_{k_1,s_2}\cdots z_{k_1,s_n} \\
&=& z_{k+l,st}z_{k_1,s_1}\cdots z_{k_1,s_n}+z_{k,s}z_{k_1+l,s_1t}z_{k_1,s_2}\cdots z_{k_1,s_n} \\
&{}& \quad +z_{k_1+l,s_1t}(z_{k,s}\ast z_{k_1,s_2}\cdots z_{k_1,s_n})+z_{k+k_1+l,ss_1t}z_{k_1,s_2}\cdots z_{k_1,s_n}.
\end{eqnarray*}
Thus, we obtain LHS$=$RHS. 

If $\mathrm{dep}(w)>1$, we can set $w=z_{k,s}u,w^{\prime}=z_{l,t}v~(\mathrm{dep}(u),\mathrm{dep}(v)>0)$. Then
\begin{eqnarray*}
\mathrm{LHS} &=& \mathcal{I}^{-1}d_{\mathcyb{sh}}\mathcal{I}(z_{k,s}u~\overline{\ast}~z_{l,t}v) \\
&=& \mathcal{I}^{-1}d_{\mathcyb{sh}}\mathcal{I}(z_{k,s}(u~\overline{\ast}~z_{l,t}v)+z_{l,t}(z_{k,s}u~\overline{\ast}~v)-z_{k+l,st}(u~\overline{\ast}~v)) \\
&=& \mathcal{I}^{-1}d_{\mathcyb{sh}}L_{z_{k,s}}M_{\frac{1}{s}}\mathcal{I}(u~\overline{\ast}~z_{l,t}v)+\mathcal{I}^{-1}d_{\mathcyb{sh}}L_{z_{l,t}}M_{\frac{1}{t}}\mathcal{I}(z_{k,s}u~\overline{\ast}~v) \\
&{}& \quad -\mathcal{I}^{-1}d_{\mathcyb{sh}}L_{z_{k+l,st}}M_{\frac{1}{st}}\mathcal{I}(u~\overline{\ast}~v) \\
&=& \mathcal{I}^{-1}L_{z_{k,s}+x^k}M_{\frac{1}{s}}d_{\mathcyb{sh}}\mathcal{I}(u~\overline{\ast}~z_{l,t}v)+\mathcal{I}^{-1}L_{z_{l,t}+x^l}M_{\frac{1}{t}}d_{\mathcyb{sh}}\mathcal{I}(z_{k,s}u~\overline{\ast}~v) \\
&{}& \quad -\mathcal{I}^{-1}L_{z_{k+l,st}+x^{k+l}}M_{\frac{1}{st}}d_{\mathcyb{sh}}\mathcal{I}(u~\overline{\ast}~v) \\
&=& L_{z_{k,s}}\mathcal{I}^{-1}d_{\mathcyb{sh}}\mathcal{I}(u~\overline{\ast}~z_{l,t}v)+L_{x^l}M_s\mathcal{I}^{-1}d_{\mathcyb{sh}}\mathcal{I}(u~\overline{\ast}~z_{l,t}v) \\
&{}& \quad +L_{z_{l,t}}\mathcal{I}^{-1}d_{\mathcyb{sh}}\mathcal{I}(z_{k,s}u~\overline{\ast}~v)+L_{x^k}M_t\mathcal{I}^{-1}d_{\mathcyb{sh}}\mathcal{I}(z_{k,s}u~\overline{\ast}~v) \\
&{}& \quad -L_{z_{k+l,st}}\mathcal{I}^{-1}d_{\mathcyb{sh}}\mathcal{I}(u~\overline{\ast}~v)-L_{x^{k+l}}M_{st}\mathcal{I}^{-1}d_{\mathcyb{sh}}\mathcal{I}(u~\overline{\ast}~v) \\
&=& L_{z_{k,s}}(\mathcal{I}^{-1}d_{\mathcyb{sh}}\mathcal{I}(u)\ast(L_{z_{l,t}}\mathcal{I}^{-1}d_{\mathcyb{sh}}\mathcal{I}(v)+L_{x^l}M_t\mathcal{I}^{-1}d_{\mathcyb{sh}}\mathcal{I}(v))) \\
&{}& \quad +L_{x^k}M_s(\mathcal{I}^{-1}d_{\mathcyb{sh}}\mathcal{I}(u)\ast(L_{z_{l,t}}\mathcal{I}^{-1}d_{\mathcyb{sh}}\mathcal{I}(v)+L_{x^l}M_t\mathcal{I}^{-1}d_{\mathcyb{sh}}\mathcal{I}(v))) \\
&{}& \quad +L_{z_{l,t}}((L_{z_{k,s}}\mathcal{I}^{-1}d_{\mathcyb{sh}}\mathcal{I}(u)+L_{x^k}M_s\mathcal{I}^{-1}d_{\mathcyb{sh}}\mathcal{I}(u))\ast\mathcal{I}^{-1}d_{\mathcyb{sh}}\mathcal{I}(v)) \\
&{}& \quad +L_{x^l}M_t((L_{z_{k,s}}\mathcal{I}^{-1}d_{\mathcyb{sh}}\mathcal{I}(u)+L_{x^k}M_s\mathcal{I}^{-1}d_{\mathcyb{sh}}\mathcal{I}(u))\ast\mathcal{I}^{-1}d_{\mathcyb{sh}}\mathcal{I}(v)) \\
&{}& \quad -L_{z_{k+l,st}}(\mathcal{I}^{-1}d_{\mathcyb{sh}}\mathcal{I}(u)\ast\mathcal{I}^{-1}d_{\mathcyb{sh}}\mathcal{I}(v)) \\
&{}& \quad -L_{x^{k+l}}M_{st}(\mathcal{I}^{-1}d_{\mathcyb{sh}}\mathcal{I}(u)\ast\mathcal{I}^{-1}d_{\mathcyb{sh}}\mathcal{I}(v)),
\end{eqnarray*}
\begin{eqnarray*}
\mathrm{RHS} &=& ((L_{z_{k,s}}+L_{x^k}M_s)\mathcal{I}^{-1}d_{\mathcyb{sh}}\mathcal{I}(u))\ast((L_{z_{l,t}}+L_{x^l}M_t)\mathcal{I}^{-1}d_{\mathcyb{sh}}\mathcal{I}(v)) \\
&=& L_{z_{k,s}}(\mathcal{I}^{-1}d_{\mathcyb{sh}}\mathcal{I}(u)\ast L_{z_{l,t}}\mathcal{I}^{-1}d_{\mathcyb{sh}}\mathcal{I}(v)) \\
&{}& \quad +L_{z_{l,t}}(L_{z_{k,s}}\mathcal{I}^{-1}d_{\mathcyb{sh}}\mathcal{I}(u)\ast \mathcal{I}^{-1}d_{\mathcyb{sh}}\mathcal{I}(v)) \\
&{}& \quad +L_{z_{k+l,st}}(\mathcal{I}^{-1}d_{\mathcyb{sh}}\mathcal{I}(u)\ast \mathcal{I}^{-1}d_{\mathcyb{sh}}\mathcal{I}(v)) \\
&{}& \quad +L_{x^k}M_s\mathcal{I}^{-1}d_{\mathcyb{sh}}\mathcal{I}(u)\ast L_{z_{l,t}}\mathcal{I}^{-1}d_{\mathcyb{sh}}\mathcal{I}(v) \\
&{}& \quad +L_{z_{k,f}}\mathcal{I}^{-1}d_{\mathcyb{sh}}\mathcal{I}(u)\ast L_{x^l}M_t\mathcal{I}^{-1}d_{\mathcyb{sh}}\mathcal{I}(v) \\
&{}& \quad +L_{x^k}M_s\mathcal{I}^{-1}d_{\mathcyb{sh}}\mathcal{I}(u)\ast L_{x^l}M_t\mathcal{I}^{-1}d_{\mathcyb{sh}}\mathcal{I}(v).
\end{eqnarray*}

Hence, we must show that 
\begin{eqnarray*}
&{}& L_{x^k}M_s(V\ast (L_{z_{l,t}}+L_{x^l}M_t)(W))+L_{z_{k,s}}(V\ast L_{x^l}M_t(W)) \\
&{}& \quad +L_{x^l}M_t((L_{z_{k,s}}+L_{x^k}M_s)(V)\ast W)+L_{z_{l,t}}(L_{x^k}M_s(V)\ast W) \\
&{}& \quad -(L_{z_{k+l,st}}+L_{x^{k+l}}M_{st})(V\ast W) \\
&{}& =L_{z_{k,s}}(V)\ast L_{x^l}M_t(W)+L_{x^k}M_s(V)\ast L_{z_{l,t}}(W) \\
&{}& \quad +L_{x^k}M_s(V)\ast L_{x^l}M_t(W)+L_{z_{k+l,st}}(V\ast W).
\end{eqnarray*}
Setting $V=z_{p,i}v, W=z_{q,j}w$, 
\begin{eqnarray*}
\mathrm{LHS} &=& L_{x^k}M_s\bigl(L_{z_{p,i}}(v\ast L_{z_{l,t}}(W))+L_{z_{l,t}}(V\ast W)+L_{z_{p+l,ti}}(v\ast W) \\
&{}& \quad +L_{z_{p,i}}(v\ast L_{x^l}M_t(W))+L_{z_{l+q,tj}}(V\ast w)+L_{z_{l+p+q,tij}}(v\ast w)\bigl) \\
&{}& \quad +L_{x^l}M_t\bigl(L_{z_{k,s}}(V\ast W)+L_{z_{q,j}}(L_{z_{k,s}}(V)\ast w)+L_{z_{k+q,sj}}(V\ast w) \\
&{}& \quad +L_{z_{k+p,si}}(v\ast W)+L_{z_{q,j}}(L_{x^l}M_s(V)\ast w)+L_{z_{k+p+q,sij}}(v\ast w)\bigl) \\
&{}& \quad +L_{z_{k,s}}(V\ast L_{x^l}M_t(W))+L_{z_{l,t}}(L_{x^k}M_s(V)\ast W) \\
&{}& \quad -L_{z_{k+l,st}}(V\ast W)-L_{x^{k+l}}M_{st}\bigl(L_{z_{p,i}}(v\ast W) \\
&{}& \quad +L_{z_{q,j}}(V\ast w)+L_{z_{p+q,ij}}(v\ast w)\bigl)
\end{eqnarray*}
and
\begin{eqnarray*}
\mathrm{RHS} &=& L_{z_{k,s}}(V\ast L_{x^l}M_t(W))+L_{z_{l+q,tj}}(L_{z_{k,s}}(V)\ast W) \\
&{}& \quad +L_{z_{k+l+q,stj}}(V\ast w)+L_{z_{k+p,si}}(v\ast L_{z_{l,t}}(W)) \\
&{}& \quad +L_{z_{l,t}}(L_{x^k}M_s(V)\ast W)+L_{z_{k+l+p,sti}}(v\ast W) \\
&{}& \quad +L_{z_{k+p,si}}(v\ast L_{x^l}M_t(W))+L_{z_{l+q,tj}}(L_{x^k}M_s(V)\ast w) \\
&{}& \quad +L_{z_{k+l+p+q,stij}}(v\ast w)+L_{z_{k+l,st}}(V\ast W).
\end{eqnarray*}
Thus, we obtain LHS$=$RHS. 

(ii)~Again, by the induction on the depth of a word. If $\mathrm{dep}(w)=1$, set $w=z_{k,s},w^{\prime}=z_{l,t}v$. If $\mathrm{dep}(v)=0$, then both sides become $z_{k+l,st}$, and the assertion holds. 

If $\mathrm{dep}(v)>0$, then
$$ \mathrm{LHS}=\mathcal{I}^{-1}d_{\mathcyb{sh}}\mathcal{I}(z_{k+l,st}v)=(L_{z_{k+l,st}}+L_{x^{k+l}}M_{st})\mathcal{I}^{-1}d_{\mathcyb{sh}}\mathcal{I}(v), $$
\begin{eqnarray*}
\mathrm{RHS} &=& z_{k,s}~\dot{\ast}~(L_{z_{l,t}}+L_{x^l}M_t)\mathcal{I}^{-1}d_{\mathcyb{sh}}\mathcal{I}(v) \\
&=& L_{z_{k+l,st}}\mathcal{I}^{-1}d_{\mathcyb{sh}}\mathcal{I}(v)+z_{k,s}~\dot{\ast}~L_{x^l}M_t\mathcal{I}^{-1}d_{\mathcyb{sh}}\mathcal{I}(v).
\end{eqnarray*}
Hence, we must show that 
$$ z_{k,s}~\dot{\ast}~L_{x^l}M_t(V)=L_{x^{k+l}}M_{st}(V) $$
Setting $V=z_{p,i}v$, we have
$$ \mathrm{LHS}=z_{k+l+p,sti}v=\mathrm{RHS}. $$
Hence, the assertion is proven. 

If $\mathrm{dep}(w)>1$, we can set $w=z_{k,s}u,w^{\prime}=z_{l,t}v~(\mathrm{dep}(u),\mathrm{dep}(v)>0)$. Then
\begin{eqnarray*}
\mathrm{LHS} &=& \mathcal{I}^{-1}d_{\mathcyb{sh}}\mathcal{I}(z_{k+l,st}(u~\overline{\ast}~v)) \\
&=& (L_{z_{k+l,st}}+L_{x^{k+l}}M_{st})\mathcal{I}^{-1}d_{\mathcyb{sh}}\mathcal{I}(u~\overline{\ast}~v) \\
&=& (L_{z_{k+l,st}}+L_{x^{k+l}}M_{st})(\mathcal{I}^{-1}d_{\mathcyb{sh}}\mathcal{I}(u)\ast \mathcal{I}^{-1}d_{\mathcyb{sh}}\mathcal{I}(v)),
\end{eqnarray*}
\begin{eqnarray*}
\mathrm{RHS} &=& (L_{z_{k,s}}+L_{x^k}M_s)\mathcal{I}^{-1}d_{\mathcyb{sh}}\mathcal{I}(u)~\dot{\ast}~(L_{z_{l,t}}+L_{x^l}M_t)\mathcal{I}^{-1}d_{\mathcyb{sh}}\mathcal{I}(v) \\
&=& L_{z_{k+l,st}}(\mathcal{I}^{-1}d_{\mathcyb{sh}}\mathcal{I}(u)\ast \mathcal{I}^{-1}d_{\mathcyb{sh}}\mathcal{I}(v)) \\
&{}& \quad +L_{z_{k,s}}\mathcal{I}^{-1}d_{\mathcyb{sh}}\mathcal{I}(u)~\dot{\ast}~L_{x^l}M_t\mathcal{I}^{-1}d_{\mathcyb{sh}}\mathcal{I}(v) \\
&{}& \quad +L_{x^k}M_s\mathcal{I}^{-1}d_{\mathcyb{sh}}\mathcal{I}(u)~\dot{\ast}~L_{z_{l,t}}\mathcal{I}^{-1}d_{\mathcyb{sh}}\mathcal{I}(v) \\
&{}& \quad +L_{x^k}M_s\mathcal{I}^{-1}d_{\mathcyb{sh}}\mathcal{I}(u)~\dot{\ast}~L_{x^l}M_t\mathcal{I}^{-1}d_{\mathcyb{sh}}\mathcal{I}(v). 
\end{eqnarray*}
Hence, we must show that
\begin{eqnarray*}
&{}& L_{z_{k,s}}(V)~\dot{\ast}~L_{x^l}M_t(W)+L_{x^k}M_s(V)~\dot{\ast}~L_{z_{l,t}}(W) \\
&{}& \quad +L_{x^k}M_s(V)~\dot{\ast}~L_{x^l}M_t(W)=L_{x^{k+l}}M_{st}(V\ast W)
\end{eqnarray*}
Setting $V=z_{p,i}v, W=z_{q,j}w$, we have
\begin{eqnarray*}
\mathrm{LHS} &=& z_{k+l+q,stj}(V\ast w)+z_{k+l+p,sti}(v\ast W)+z_{k+l+p+q,stij}(v\ast w) \\
&=& L_{x^{k+l}}M_{st}\bigl(z_{q,j}(V\ast w)+z_{p,i}(v\ast W)+z_{p+q,ij}(v\ast w)\bigl) \\
&=& \mathrm{RHS}.
\end{eqnarray*}
This completes the proof. 
\end{proof}

\begin{proof}[Proof of Lemma \ref{lem3-8}]
First, we show the identity
\begin{equation}
F_sL_{z_{k,t}}=L_{z_{k,\frac{1-\delta(st)st}{1-\delta(s)s}}}F_{st} \label{eq16}
\end{equation}
for $s,t\in\Lambda$. If $s=1$, then \\
$ \displaystyle\quad F_1L_{z_{k,t}}=\mathcal{I}^{-1}\iota\mathcal{I}L_{z_{k,t}}=\mathcal{I}^{-1}\iota L_{z_{k,t}}\mathcal{I}M_t=\mathcal{I}^{-1}(\delta(t)L_{z_{k,1-t}}+(1-\delta(t))L_{z_{k,1}})\iota \mathcal{I}M_t $ \\
$ \displaystyle\qquad =\delta(t)L_{z_{k,1-t}}M_{\frac{1}{1-t}}\mathcal{I}^{-1}\iota \mathcal{I}M_t+(1-\delta(t))L_{z_{k,1}}\mathcal{I}^{-1}\iota \mathcal{I}M_t=L_{z_{k,1-\delta(t)t}}F_t. $ \\
If $s\neq 1$, then \\
$\displaystyle\quad F_sL_{z_{k,t}}=M_{\frac{1}{1-s}}\mathcal{I}^{-1}\iota\mathcal{I}M_sL_{z_{k,t}}=M_{\frac{1}{1-s}}\mathcal{I}^{-1}\iota L_{z_{k,st}}\mathcal{I}M_{st} $ \\
$\displaystyle\qquad =M_{\frac{1}{1-s}}\mathcal{I}^{-1}(\delta(st)L_{z_{k,1-st}}+(1-\delta(st))L_{z_{k,1}})\iota \mathcal{I}M_{st} $ \\
$ \displaystyle\qquad =\delta(st)L_{z_{k,\frac{1-st}{1-s}}}M_{\frac{1}{1-st}}\mathcal{I}^{-1}\iota \mathcal{I}M_{st}+(1-\delta(st))L_{z_{k,\frac{1}{1-s}}}\mathcal{I}^{-1}\iota \mathcal{I}M_{st}=L_{z_{k,\frac{1-\delta(st)st}{1-s}}}F_{st}. $ \\
Therefore, we have the assertion (\ref{eq16}). 

To prove Lemma, we must show that
\begin{equation}
F_s(w)~\overline{\ast}~w^{\prime}=F_s(w~\overline{\ast}~w^{\prime}) \label{eq17}
\end{equation}
based on Lemma \ref{lem3-7} i)${}^{\prime}$. The equality is proven by induction on the total depth of $w$ and $w^{\prime}$. It is simple to show that $\mathrm{dep}(w)=1$ or $\mathrm{dep}(w^{\prime})=1$, hence (\ref{eq17}) holds for $\mathrm{dep}(w)+\mathrm{dep}(w^{\prime})\le 1$. Set $w=z_{k_1,t}w_1$, $w^{\prime}=z_{k^{\prime}_1,1}w^{\prime}_1$ ($w\in\mathcal{A}_{\Lambda,>0}^1,w^{\prime}\in\mathcal{A}_{\{1\},>0}^1$). According to the identity (\ref{eq16}) and the induction hypothesis, we have \\
$\displaystyle\quad F_s(w)~\overline{\ast}~w^{\prime}=L_{z_{k_1,\frac{1-\delta(st)st}{1-\delta(s)s}}}F_{st}(w_1)~\overline{\ast}~L_{z_{k^{\prime}_1,1}}(w^{\prime}_1)$ \\
$\displaystyle\qquad =L_{z_{k_1,\frac{1-\delta(st)st}{1-\delta(s)s}}}(F_{st}(w_1)~\overline{\ast}~L_{z_{k^{\prime}_1,1}}w^{\prime}_1)+L_{z_{k^{\prime}_1,1}}(F_s(w)~\overline{\ast}~w^{\prime}_1)-L_{z_{k_1+k^{\prime}_1,\frac{1-\delta(st)st}{1-\delta(s)s}}}(F_{st}(w_1)~\overline{\ast}~w^{\prime}_1)$ \\
$\displaystyle\qquad =L_{z_{k_1,\frac{1-\delta(st)st}{1-\delta(s)s}}}(F_{st}(w_1~\overline{\ast}~L_{z_{k^{\prime}_1,1}}w^{\prime}_1))+L_{z_{k^{\prime}_1,1}}(F_s(w~\overline{\ast}~w^{\prime}_1))-L_{z_{k_1+k^{\prime}_1,\frac{1-\delta(st)st}{1-\delta(s)s}}}(F_{st}(w_1~\overline{\ast}~w^{\prime}_1))$ \\
$\displaystyle\qquad =F_sL_{z_{k_1,t}}(w_1~\overline{\ast}~z_{k^{\prime}_1,1}w^{\prime}_1)+F_sL_{z_{k^{\prime}_1,1}}(w~\overline{\ast}~w^{\prime}_1)-F_sL_{z_{k_1+k^{\prime}_1,t}}(w_1~\overline{\ast}~w^{\prime}_1)$ \\
$\displaystyle\qquad =F_s(w~\overline{\ast}~w^{\prime}).$ \\
Hence, (\ref{eq17}) holds. (Actually, (\ref{eq17}) holds if $\overline{\ast}$ is changed to $\ast$.)
\end{proof}

\begin{proof}[Proof of Lemma \ref{lem4-3}]
The proof is given by induction on $n$. Since $\widehat{\psi}_1^{(c)}(x)=L_xL_{y_1}$, $\widehat{\psi}_1^{(c)}(y_1)=-L_xL_{y_1}$ and $\widehat{\psi}_1^{(c)}(y_s)=-L_xL_{y_s}+L_{y_s}L_{y_1-y_s}(s\neq 1)$, we obtain the lemma for $n=1$ by setting $\widehat{\phi}_0^{(c)}=\mathrm{id}_{\mathcal{A}_{\mu_r}}$. 

Suppose that the lemma is proven for $n$. By the recursive rule of $\widehat{\psi}_n^{(c)}(u)$ and the induction hypothesis, we have
\begin{eqnarray*}
n\widehat{\psi}_{n+1}^{(c)}(u)
&=&[\widehat{\theta}^{(c)},\widehat{\psi}_n^{(c)}(u)]-\frac{1}{2}\bigl(L_z\widehat{\psi}_n^{(c)}(u)+\widehat{\psi}_n^{(c)}(u)L_z\bigl)-c\widehat{\psi}_n^{(c)}(u)\partial_1 \\
&=&[\widehat{\theta}^{(c)},(-1)^{\nu(u)}\bigl(L_x\widehat{\phi}_{n-1}^{(c)}L_{y_1+\nu(u)(u-y_1)}+\nu(u)L_u\widehat{\phi}_{n-1}^{(c)}L_{u-y_1}\bigl)] \\
\displaystyle &{}& \quad -\frac{1}{2}(-1)^{\nu(u)}L_z\bigl(L_x\widehat{\phi}_{n-1}^{(c)}L_{y_1+\nu(u)(u-y_1)}+\nu(u)L_u\widehat{\phi}_{n-1}^{(c)}L_{u-y_1}\bigl) \\
\displaystyle &{}& \quad -\frac{1}{2}(-1)^{\nu(u)}\bigl(L_x\widehat{\phi}_{n-1}^{(c)}L_{y_1+\nu(u)(u-y_1)}+\nu(u)L_u\widehat{\phi}_{n-1}^{(c)}L_{u-y_1}\bigl)L_z \\
\displaystyle &{}& \quad -c(-1)^{\nu(u)}\bigl(L_x\widehat{\phi}_{n-1}^{(c)}L_{y_1+\nu(u)(u-y_1)}+\nu(u)L_u\widehat{\phi}_{n-1}^{(c)}L_{u-y_1}\bigl)\partial_1.
\end{eqnarray*}
Write the leading Lie bracket term as the sum
\begin{eqnarray*}
\displaystyle &{}& (-1)^{\nu(u)}\Bigl\{[\widehat{\theta}^{(c)},L_x]\widehat{\phi}_{n-1}^{(c)}L_{y_1+\nu(u)(u-y_1)}+L_x[\widehat{\theta}^{(c)},\widehat{\phi}_{n-1}^{(c)}]L_{y_1+\nu(u)(u-y_1)} \\
\displaystyle &{}& \quad +L_x\widehat{\phi}_{n-1}^{(c)}[\widehat{\theta}^{(c)},L_{y_1+\nu(u)(u-y_1)}]+\nu(u)\bigl([\widehat{\theta}^{(c)},L_u]\widehat{\phi}_{n-1}^{(c)}L_{u-y_1} \\
\displaystyle &{}& \quad +L_u[\widehat{\theta}^{(c)},\widehat{\phi}_{n-1}^{(c)}]L_{u-y_1}+L_u\widehat{\phi}_{n-1}^{(c)}[\widehat{\theta}^{(c)},L_{u-y_1}]\bigl)\Bigl\}.
\end{eqnarray*}
Using the identity $[\widehat{\theta}^{(c)},L_u]=L_{\widehat{\theta}^{(c)}(u)}+cL_u\partial_1$, which can be easily shown, we have 
\begin{eqnarray*}
&{}& n\widehat{\psi}_{n+1}^{(c)}(u) \\
&=& (-1)^{\nu(u)}\Bigl\{(L_{\widehat{\theta}^{(c)}(x)}+cL_x\partial_1)\widehat{\phi}_{n-1}^{(c)}L_{y_1+\nu(u)(u-y_1)}+L_x[\widehat{\theta}^{(c)},\widehat{\phi}_{n-1}^{(c)}]L_{y_1+\nu(u)(u-y_1)} \\
\displaystyle &{}& \quad +L_x\widehat{\phi}_{n-1}^{(c)}(L_{\widehat{\theta}^{(c)}(y_1+\nu(u)(u-y_1))}+cL_{y_1+\nu(u)(u-y_1)}\partial_1) \\
\displaystyle &{}& \quad +\nu(u)\bigl((L_{\widehat{\theta}^{(c)}(u)}+cL_u\partial_1)\widehat{\phi}_{n-1}^{(c)}L_{u-y_1} \\
\displaystyle &{}& \quad +L_u[\widehat{\theta}^{(c)},\widehat{\phi}_{n-1}^{(c)}]L_{u-y_1}+L_u\widehat{\phi}_{n-1}^{(c)}(L_{\widehat{\theta}^{(c)}(u-y_1)}+cL_{u-y_1}\partial_1)\bigl) \\
\displaystyle &{}& \quad -\frac{1}{2}L_z\bigl(L_x\widehat{\phi}_{n-1}^{(c)}L_{y_1+\nu(u)(u-y_1)}+\nu(u)L_u\widehat{\phi}_{n-1}^{(c)}L_{u-y_1}\bigl) \\
\displaystyle &{}& \quad -\frac{1}{2}\bigl(L_x\widehat{\phi}_{n-1}^{(c)}L_{y_1+\nu(u)(u-y_1)}+\nu(u)L_u\widehat{\phi}_{n-1}^{(c)}L_{u-y_1}\bigl)L_z \\
\displaystyle &{}& \quad -c\bigl(L_x\widehat{\phi}_{n-1}^{(c)}L_{y_1+\nu(u)(u-y_1)}+\nu(u)L_u\widehat{\phi}_{n-1}^{(c)}L_{u-y_1}\bigl)\partial_1\Bigl\}. 
\end{eqnarray*}
Considering $\theta^{(c)}(u)=\frac{1}{2}(uz+zu)$, 
\begin{eqnarray*}
\displaystyle n\widehat{\psi}_{n+1}^{(c)}(u) &=& (-1)^{\nu(u)}\Bigl\{\frac{1}{2}L_xL_z\widehat{\phi}_{n-1}^{(c)}L_{y_1+\nu(u)(u-y_1)}+cL_x\partial_1\widehat{\phi}_{n-1}^{(c)}L_{y_1+\nu(u)(u-y_1)} \\
\displaystyle &{}& \quad +L_x[\widehat{\theta}^{(c)},\widehat{\phi}_{n-1}^{(c)}]L_{y_1+\nu(u)(u-y_1)}+\frac{1}{2}L_x\widehat{\phi}_{n-1}^{(c)}L_zL_{y_1+\nu(u)(u-y_1)} \\
\displaystyle &{}& \quad +\nu(u)\Bigl(\frac{1}{2}L_uL_z\widehat{\phi}_{n-1}^{(c)}L_{u-y_1}+cL_u\partial_1\widehat{\phi}_{n-1}^{(c)}L_{u-y_1} \\
\displaystyle &{}& \quad +L_u[\widehat{\theta}^{(c)},\widehat{\phi}_{n-1}^{(c)}]L_{u-y_1}+\frac{1}{2}L_u\widehat{\phi}_{n-1}^{(c)}L_zL_{u-y_1}\Bigl)\Bigl\}.
\end{eqnarray*}
By setting
\begin{equation}
\displaystyle \widehat{\phi}_n^{(c)}=\frac{1}{n}\Bigl([\widehat{\theta}^{(c)},\widehat{\phi}_{n-1}^{(c)}]+\frac{1}{2}\bigl(L_z\widehat{\phi}_{n-1}^{(c)}+\widehat{\phi}_{n-1}^{(c)}L_z\bigl)+c\partial_1\widehat{\phi}_{n-1}^{(c)}\Bigl), \label{eq14}
\end{equation}
we have $n\widehat{\psi}_{n+1}^{(c)}(u)=(-1)^{\nu(u)}n\bigl(L_x\widehat{\phi}_n^{(c)}L_{y_1+\nu(u)(u-y_1)}+\nu(u)L_u\widehat{\phi}_n^{(c)}L_{u-y_1}\bigl).$ This shows the lemma. 
\end{proof}


\section*{Appendix: Tables}\label{sec6}

We give the maximal number of linearly independent relations supplied by each set of relations among MLV's for each $\mu_r$ with $1\le r\le 6$. The lineup is the derivation relation (``$\partial_n$''), the extended derivation relation (``$\partial_n^{(c)}$'') for all $c\in\mathbb{Q}$, the linear part of Corollary \ref{cor3-13} (``lin.'') sequentially from the top. (We have already proven ``$\partial_n$''$\subset$``$\partial_n^{(c)}$''$\subset$``lin.'') Computations were performed using Risa/Asir, an open-source general computer algebra system. 

\begin{center}
\underline{$r=1$} (MZV case)

\begin{tabular}{|c||c|c|c|c|c|c|c|c|c|c|c|c|}
\hline
weight & 3 & 4 & 5 & 6 & 7 & 8 & 9 & 10 & 11 & 12 & 13 & 14  \\
\hhline{|=#=|=|=|=|=|=|=|=|=|=|=|=|}
$\partial_n$ & 1 & 2 & 5 & 10 & 22 & 44 & 90 & 181 & 363 & 727 & 1456 & 2912  \\
\hline
$\partial_n^{(c)}$ & 1 & 2 & 5 & 10 & 23 & 46 & 98 & 200 & 410 & 830 & 1679 & $\cdots$  \\
\hline
lin. & 1 & 2 & 5 & 10 & 23 & 46 & 98 & 200 & 413 & 838 & 1713 & $\cdots$   \\
\hline
$\#$\{index set\} & 2 & 4 & 8 & 16 & 32 & 64 & 128 & 256 & 512 & 1024 & 2048 & 4096  \\
\hline
\end{tabular}
\end{center}
\begin{center}
\underline{$r=2$}

\begin{tabular}{|c||c|c|c|c|c|c|}
\hline
weight & 3 & 4 & 5 & 6 & 7 & 8   \\
\hhline{|=#=|=|=|=|=|=|}
$\partial_n$ & 4 & 14 & 46 & 140 & 426 & 1280   \\
\hline
$\partial_n^{(c)}$ & 4 & 14 & 48 & 150 & 464 & 1402   \\
\hline
lin. & 4 & 14 & 48 & 150 & 468 & 1422   \\
\hline
$\#$\{index set\} & 12 & 36 & 108 & 324 & 972 & 2916   \\
\hline
\end{tabular}
\end{center}
\begin{center}
\underline{$r=3$}

\begin{tabular}{|c||c|c|c|c|}
\hline
weight & 3 & 4 & 5 & 6    \\
\hhline{|=#=|=|=|=|}
$\partial_n$ & 9 & 42 & 177 & 714   \\
\hline
$\partial_n^{(c)}$ & 9 & 42 & 183 & 750   \\
\hline
lin. & 9 & 42 & 183 & 750  \\
\hline
$\#$\{index set\} & 36 & 144 & 576 & 2304  \\
\hline
\end{tabular}
\end{center}
\begin{center}
\underline{$r=4$}

\begin{tabular}{|c||c|c|c|}
\hline
weight & 3 & 4 & 5    \\
\hhline{|=#=|=|=|}
$\partial_n$ & 16 & 92 & 476  \\
\hline
$\partial_n^{(c)}$ & 16 & 92 & 488  \\
\hline
lin. & 16 & 92 & 488  \\
\hline
$\#$\{index set\} & 80 & 400 & 2000  \\
\hline
\end{tabular}
\end{center}
\begin{center}
\underline{$r=5$}

\begin{tabular}{|c||c|c|c|}
\hline
weight & 3 & 4 & 5   \\
\hhline{|=#=|=|=|}
$\partial_n$ & 25 & 170 & 1045  \\
\hline
$\partial_n^{(c)}$ & 25 & 170 & 1065   \\
\hline
lin. & 25 & 170 & 1065   \\
\hline
$\#$\{index set\} & 150 & 900 & 5400   \\
\hline
\end{tabular}
\end{center}
\begin{center}
\underline{$r=6$}

\begin{tabular}{|c||c|c|}
\hline
weight & 3 & 4   \\
\hhline{|=#=|=|}
$\partial_n$ & 36 & 282   \\
\hline
$\partial_n^{(c)}$ & 36 & 282   \\
\hline
lin. & 36 & 282   \\
\hline
$\#$\{index set\} & 252 & 1764   \\
\hline
\end{tabular}
\end{center}



\end{document}